\documentclass[12pt,reqno]{amsart}
\usepackage{amsmath}
\usepackage{amssymb}
\usepackage{verbatim}
\usepackage{mathrsfs}
\usepackage{graphicx} %If you want to include postscript graphics
\usepackage{caption}
\usepackage{subcaption}
\usepackage[noadjust]{cite}
\usepackage{hyperref}
\usepackage{color}

\usepackage[noadjust]{cite}

\usepackage[top=0.5in,bottom=0.6in,left=0.67in,right=0.67in]{geometry}

\newtheorem{lemma}{Lemma}

\newtheorem{theorem}[lemma]{Theorem}
\newtheorem{remark}[lemma]{Remark}

\newtheorem{exmp}{Example}

\newtheorem{definition}{Definition}

\numberwithin{equation}{section}
\newcommand{\C}{\mathbb{C}}    %complex number field
\newcommand{\N}{\mathbb{N}}    %natural numbers
\newcommand{\R}{\mathbb{R}}    %real number field
\newcommand{\Z}{\mathbb{Z}}    %integers

\newcommand{\wh}{\widehat}
\renewcommand{\le}{\leqslant}
\renewcommand{\ge}{\geqslant}

\newcommand{\ol}{\overline}

\newcommand{\bo}{\mathscr{O}} %standard big O notation

\newcommand{\eps}{\epsilon}

\newcommand{\er}{\eqref}

\newcommand{\bp}{ \begin{proof} }
	\newcommand{\ep}{\hfill \end{proof} }
\newcommand{\be}{ \begin{equation} }
\newcommand{\ee}{ \end{equation} }
\newcommand{\tp}{\mathsf{T}}
\newcommand{\dm}{\mathsf{M}} %or M
\newcommand{\dn}{\mathsf{N}} %or N
\newcommand{\sr}{\operatorname{sr}}  %sum rules
\newcommand{\sm}{\operatorname{sm}}  %smoothness exponent
\newcommand{\lp}[1]{l_{#1}(\mathbb{Z})}

\newcommand{\td}{\boldsymbol{\delta}}  %Dirac/Kronicker sequence

\newcommand{\sd}{\mathcal{S}}  %subdivision operator S
\newcommand{\dR}{\mathbb{R}^d}

\newcommand{\dZ}{\mathbb{Z}^d}

\newcommand{\dLp}[1]{L_{#1}(\mathbb{R}^d)}

\newcommand{\fsupp}{\text{fsupp}}

\newcommand{\vertiii}[1]{{\left\vert\kern-0.25ex\left\vert\kern-0.25ex\left\vert #1
		 \right\vert\kern-0.25ex\right\vert\kern-0.25ex\right\vert}}

\newcommand{\dlp}[1]{l_{#1}(\mathbb{Z}^d)}

\newcommand{\dlpp}[1]{l_{#1}(\mathbb{Z}^2)}

\begin{document}
	
%\title[Balanced Multivariate Quasi-tight Multiframelets] {Balanced Multivariate Quasi-tight Multiframelets Derived From Arbitrary Refinable Vector Functions}

\title[Quasi-stationary Subdivision Schemes] {Quasi-stationary Subdivision Schemes in Arbitrary Dimensions}

\author{Ran Lu}
\address{College of Science, Hohai University, Nanjing, China 211100.
	\quad  {\tt rlu3@hhu.edu.cn}}

\author{Bin Han}
\address{Department of Mathematical and Statistical Sciences,
	University of Alberta, Edmonton,\quad Alberta, Canada T6G 2G1.
	\quad {\tt bhan@ualberta.ca} }

\thanks{The research of the first author was supported by the National Natural Science Foundation of China under grant 12201178. }

\thanks{The research of the second author was supported by the Natural Sciences and Engineering Research Council of Canada (NSERC) under grant RGPIN 2024-04991.
}	
	
	%\thanks{Contact information of corresponding author Bin Han: E-mail: bhan@ualberta.ca, Phone: 1-587-8828375, Fax: 1-780-4926826,  Web: http://www.ualberta.ca/$\sim$bhan}
	
	\makeatletter \@addtoreset{equation}{section} \makeatother
	
	\begin{abstract}Stationary subdivision schemes have been extensively studied and have numerous applications in CAGD and wavelet analysis. To have high-order smoothness of the scheme, it is usually inevitable to enlarge the support of the mask that is used, which is a major difficulty with stationary subdivision schemes due to complicated implementation and dramatically increased special subdivision rules at extraordinary vertices. In this paper, we introduce the notion of a multivariate quasi-stationary subdivision scheme and fully characterize its convergence and smoothness. We will also discuss the general procedure of designing interpolatory masks with short support that yields smooth quasi-stationary subdivision schemes. Specifically, using the dyadic dilation of both triangular and quadrilateral meshes, for each smoothness exponent $m=1,2$, we obtain examples of $C^m$-convergent quasi-stationary $2I_2$-subdivision schemes with bivariate symmetric masks having at most $m$-ring stencils. Our examples demonstrate the advantage of quasi-stationary subdivision schemes, which can circumvent the difficulty above with stationary subdivision schemes.
\end{abstract}


\maketitle

\pagenumbering{arabic}

\section{Introduction and Motivation}

Subdivision schemes are fast iterative averaging algorithms for computing refinable functions and wavelets. Due to the cascade/multi-scale structure and intrinsic connections to splines and wavelets, subdivision schemes are of great interest in many applications such as computer-aided geometric design (CAGD) for generating smooth curves and surfaces (\cite{dd89,dl02,dlg87,dlg90,fhs22,hj98,hj06,hyx05,ms19}), solving PDEs numerically using multi-scale methods (\cite{hm21} and references therein), and data processing with discrete wavelet/framelet transforms (\cite{hanbook,hl19pp}). This paper focuses on \emph{multivariate quasi-stationary subdivision schemes} in arbitrary dimensions. In this section, we recall some basics of subdivision schemes and explain the motivations and contributions of our work.

\subsection{Stationary subdivision schemes: convergence and smoothness}

The classical way to implement a subdivision scheme is by performing a subdivision operation using the same mask at every level, and such a scheme is often known as a stationary subdivision scheme. To be specific, let us introduce some basic notations. By a \emph{$d$-dimensional mask/filter} we mean a sequence $a=\{a(k)\}_{k\in\dZ}:\dZ\to \C$ such that $a(k)\ne 0$ for only finitely many terms. By $\dlp{0}$ we denote the linear space of all $d$-dimensional filters. To perform a subdivision scheme, we need a mask $a\in\dlp{0}$ normalized by
\be\label{wha:1}
\sum_{k\in\dZ}a(k)=1,
\ee
and we shall use a \emph{dilation matrix} $\dm I_d$ where $\dm\in\N\setminus\{1\}$ and $I_d$ denotes the $d\times d$ identity matrix. The \emph{$\dm I_d$-subdivision operator} $\sd_{a,\dm I_d}$ that uses the mask $a$ is then defined by
$$
[\sd_{a,\dm I_d}v](k):=\dm^d\sum_{n\in\dZ}v(n)a(k-\dm n),\quad\forall v\in\dlp{0},\,k\in\dZ.
$$
Given an initial data $v\in\dlp{0}$, the  \emph{stationary $\dm I_d$-subdivision scheme} that employs the mask $a$ iteratively generates a sequence $\{\sd_{a,\dm I_d}^nv\}_{n=1}^\infty$, which is expected to converge, after locating the value $[\sd_{a,\dm I_d}^nv](k)$ at the position $\dm^{-n}k$ for $k\in \dZ$, to a smooth function/surface $\eta_v$.

Stationary subdivision schemes are closely related to \emph{refinable functions}. For a mask $a\in\dlp{0}$ satisfying \er{wha:1} and for a dilation matrix $\dm I_d$, it is well known that there exists a compactly supported distribution $\phi$ such that the following \emph{$\dm I_d$-refinement equation} holds:
\be\label{ref}
\phi(x)=\dm^d \sum_{k\in\dZ}a(k)\phi(\dm x-k),\qquad\forall x\in\dR.
\ee
Any $\phi$ satisfying \er{ref} is called an \emph{$\dm I_d$-refinable function of the mask $a$}. In such cases, the mask $a$ is called the \emph{$\dm I_d$-refinement mask of the distribution $\phi$}. One powerful tool for studying refinable properties is the Fourier transform. For $f\in\dLp{1}$, its \emph{Fourier transform} is defined by
\[
\wh{f}(\xi):=\int_{\dR}f(x)e^{-ix\cdot\xi}dx,\qquad\forall\xi\in\dR.
\]
The definition of the Fourier transform is naturally extended to tempered distributions. For any $a\in\dlp{0}$, define its \emph{Fourier series} by
\[
\wh{a}(\xi):=\sum_{k\in\dZ}a(k)e^{-ik\cdot\xi},\qquad \forall\xi\in\dR.
\]
Using the Fourier transform, \er{ref} is equivalent to
\be\label{ref:f}
\wh{\phi}(\dm\xi)=\wh{a}(\xi)\wh{\phi}(\xi),\qquad\forall\xi\in\dR.
\ee
Suppose $a\in\dlp{0}$ satisfies $\wh{a}(0)=1$ (this condition is equivalent to $\sum_{k\in \dZ} a(k)=1$ in \er{wha:1}), one can define a compactly supported distribution $\phi$ through its Fourier transform by
\be\label{ref:a}\wh{\phi}(\xi):=\prod_{j=1}^\infty \wh{a}(\dm^{-j}\xi),\qquad\forall\xi\in\dR.
\ee
Then clearly $\phi$ satisfies \er{ref} with $\wh{\phi}(0)=1$. In this case, the above function $\phi$ is called the \emph{standard $\dm I_d$-refinable function} of the mask $a$. Generally, an $\dm I_d$-refinable function $\phi$ of a mask $a\in\dlp{0}$ does not have an analytic explicit expression. Fortunately, one can approximate $\phi$ by implementing the $\dm I_d$-subdivision scheme with its refinement mask $a$, provided that the scheme is \emph{convergent}. Let $a\in\dlp{0}$ be such that $\wh{a}(0)=1$ and $\dm\in \N\setminus\{1\}$. If for every $v\in\dlp{0}$, there exists a continuous $d$-dimensional function $\eta_v$ on $\dR$ such that
\[
\lim_{n\to\infty}\sup_{k\in\dZ}|[\sd_{a,\dm I_d}^nv](k)-\eta_v(\dm^{-n}k)|=0,
\]
then we say that the $\dm I_d$-subdivision scheme with the $d$-dimensional mask $a$ is \emph{convergent}. It is known (e.g., see \cite[Theorem~4.3]{han03} or \cite[Theorem~7.3.1]{hanbook}) that the $\dm I_d$-subdivision scheme with a mask $a\in \dlp{0}$ is convergent if and only if $\sm_\infty(a,\dm I_d)>0$ (see \eqref{sm} for its definition).
For a convergent $\dm I_d$-subdivision scheme employing a mask $a\in\dlp{0}$, if $v(k)=\td(k)$ for all $k\in\dZ$, where
\be\label{td}\td(x):=\begin{cases}1, &x=0,\\
	0, & x\in\dR\setminus\{0\},\end{cases}
\ee
then the limit function $\eta_{\td}$ is the standard $\dm I_d$-refinable function $\phi$ that is defined via \er{ref:a}.

For a convergent subdivision scheme, if any initial data $v\in\dlp{0}$ can be interpolated by its limit function $\eta_v$, that is,
\[
\eta_v(k)=v(k),\qquad\forall k\in\dZ,
\]
then the subdivision scheme is \emph{$\dm I_d$-interpolatory}. Interpolatory subdivision schemes are important in sampling theory, signal processing, and CAGD. The literature has extensively studied their theoretical properties and applications (e.g., see \cite{dlg87,dlg90,cj05,cj06,dhmm05,han23-2,hj98b,hj06,jlz09,jo03,grv22} and many references therein). For a convergent $\dm I_d$-subdivision scheme that employs the mask $a\in\dlp{0}$,  using the linearity of the subdivision operator, the limit function $\eta_v$ of any input data $v\in\dlp{0}$ must satisfy
\[
\eta_v(x)=\sum_{k\in\dZ}v(k)\phi(x-k),\qquad\forall x\in\dR,
\]
where $\phi$ is the standard $\dm I_d$-refinable function associated with the mask $a$ that is defined via \er{ref:a}. Therefore, a convergent $\dm I_d$-subdivision scheme is $\dm I_d$-interpolatory if and only if the $\dm I_d$-refinable function $\phi$ is interpolatory, that is, $\phi$ is continuous and
\be\label{int:phi}
\phi(k)=\td(k),\quad\forall k\in\Z^d.
\ee
Furthermore, for a convergent subdivision scheme, \er{int:phi} implies
\be\label{int:a}
a(\dm k)=\dm^{-d}\td(k),\quad \forall k\in\Z^d,
\ee
that is, $a$ must be an \emph{$\dm I_d$-interpolatory mask}.
For an $\dm I_d$-refinable function $\phi$ with a finitely supported mask $a\in \dlp{0}$, it is known (e.g., see \cite[Corollary~5.2]{han03}) that $\phi$ is interpolating if and only if $\sm_\infty(a,\dm I_d)>0$ and the mask $a$ is $\dm I_d$-interpolatory.

In applications such as CAGD, to have good visual quality of the generated subdivision curves or surfaces, a stationary subdivision scheme of high-order smoothness is desired. Smooth stationary subdivision schemes and their applications have been well-studied in the literature; see, for instance, \cite{cdm91,cjd02,dl02,hanbook,hj98,hj06} and many references therein. To define the smoothness of a stationary subdivision scheme, for $h\in\dZ$, we first define the \emph{backward difference operator} $\nabla_h$ by:
\[
\nabla_h v(k)=v(k)-v(k-h),\quad\forall v\in \dlp{0},\,k\in\dZ.
\]
For every $\mu=(\mu_1,\dots,\mu_d)\in\N_0^d$, we define
\[
\nabla^\mu:=\nabla^{\mu_1}_{e_1}\dots\nabla^{\mu_d}_{e_d},
\]
where $e_j$ is the $j$-th coordinate vector in $\dR$ for all $j=1,\dots,d$. For any $\mu\in\N_0^d$ and $v\in\dlp{0}$, observe that $\nabla^\mu v =[\nabla^\mu \td]*v$, where the convolution of two filters $v,u\in\dlp{0}$ is defined as
\[
[v*u](k):=\sum_{q\in\dZ}v(k-q)u(q),\qquad\forall k\in\dZ.
\]
It is straightforward to check that
\[
\wh{\nabla^\mu v}(\xi)=\wh{\nabla^\mu\td}(\xi)\wh{v}(\xi)=(1-e^{-i\xi_1})^{\mu_1}
\dots(1-e^{-i\xi_d})^{\mu_d}\wh{v}(\xi),
\quad\forall\xi=(\xi_1,\dots,\xi_d)\in\dR.
\]
Next, we recall the notion of \emph{$C^m$-convergence} where $m\in\N_0$ measures the order of smoothness. Let $s\in\N_0$ and $a\in\dlp{0}$ satisfying $\wh{a}(0)=1$. If for every initial data $v\in\dlp{0}$, there exists $\eta_v\in C^m(\dR)$ such that
\be\label{cm:conv}\lim_{n\to\infty}\sup_{k\in\dZ}\left|\dm^{jn}[\nabla^\mu\sd_{a,\dm I_d}^nv](k)-\partial^\mu\eta_v(\dm^{-n}k)\right|=0,\quad\forall  \mu\in\N_{0,j}^d,\,j=0,1,\dots,m, \ee
where $\N_{0,j}^d:=\{\mu\in\N_0^d:|\mu|=j\},$ then we say that the $\dm I_d$-subdivision scheme that employs the mask $a$ is \emph{$C^m$-convergent}. The smoothness of a stationary subdivision scheme is fully characterized by its underlying mask $a$. To do this, we need to introduce two technical quantities: the \emph{sum rule orders} and the \emph{smoothness exponenets} of the mask $a$. Let $f$ and $g$ be smooth functions and $m\in\N_0$, recall the following \emph{big $\bo$ notation}:
\[
f(\xi)=g(\xi)+\bo(\|\xi-\xi_0\|^m),\quad \xi\to\xi_0,
\]
which means
$$\partial^\mu f(\xi_0)=\partial^\mu g(\xi_0),\quad\forall \mu\in\N_0^d\text{ such that }|\mu|<m.$$
Let $a\in\dlp{0}$ be a $d$-dimensional filter and $\dm\in\N\setminus\{1\}$.
\begin{enumerate}
	\item[(1)] For $m\in\N_0$, we say that the mask $a$ has \emph{order $m$ sum rules} with respect to $\dm I_d$ if
	\be\label{sr:def}\sum_{k\in\dZ}(\gamma+\dm k)^\mu a(\gamma+\dm k)=\dm^{-d}\sum_{k\in\dZ}k^\mu a(k),\quad\forall \gamma\in\dZ,\,\mu\in\N_{0}^d\text{ such that }|\mu|<m,\ee	
	or equivalently, 	
	 \be\label{sr}\wh{a}(\xi+2\pi\omega)=\bo(\|\xi\|^m),\quad \xi\to 0,\qquad\forall\omega\in\Omega_{\dm}\setminus\{0\},\ee
	where
	 \be\label{om}\Omega_{\dm}:=[\dm^{-1}\dZ]\cap[0,1)^d.\ee
	Define
	\[
	\sr(a,\dm I_d):=\sup\{m\in\N_0:\,\text{\er{sr:def} or \er{sr} holds}\}.
	\]
	
	\item[(2)] Let $1\le p\le\infty$ and suppose $\sr(a,\dm I_d)=m$. We define
	%\emph{$L_p$-joint spectral radius} of $a$ with respect to $\dm I_d$ by
	\[
	\rho_m(a,\dm I_d)_p:=\sup\left\{\limsup_{n\to\infty}\left\|\nabla^\mu\sd_{a,\dm}^n\td\right\|_{l_p(\dZ)}^{1/n}:\,\mu\in\N_{0,m}^d\right\},
	\]
	and the \emph{$L_p$-smoothness exponent} of the mask $a$ (with respect to $\dm I_d$) by
\be \label{sm}
\sm_p(a,\dm I_d):=\frac{d}{p}-\log_\dm[\rho_m(a,\dm I_d)_p].
\ee
\end{enumerate}
By \cite[Theorem~4.3]{han03} (also see
\cite[Theorem 2.1]{hj06} and \cite[Theorem~7.3.1]{hanbook}), the stationary $\dm I_d$-subdivision scheme that employs the mask $a$ is $C^m$-convergent if and only if $\sm_\infty(a,\dm I_d)>m$. In particular, for a $C^m$-convergent stationary subdivision scheme, the standard $\dm I_d$-refinable function $\phi$ derived from the mask $a$ via \er{ref:a} belongs to $C^m(\dR)$. Furthermore,the partial derivative $\partial^\mu\phi$ is compactly supported and uniformly continuous for all $\mu\in\N_0^d$ with $|\mu|\le m$. Generally, there is no efficient way to compute $\sm_\infty(a,\dm I_d)$. One way to estimate $\sm_\infty(a,\dm I_d)$ is from $\sm_2(a,\dm I_d)$ which can be efficiently computed (e.g. see \cite{han03,han03-2,jj03}). Using the definition of the $L_p$-smoothness exponent, one can directly obtain the following lower bound of $\sm_\infty(a,\dm I_d)$ (e.g., see \cite[Theorem 3.1]{han03-2} or \cite[Lemma 3.1]{lu24}):
\be\label{sm:inf:2}\sm_\infty(a,\dm I_d)\ge \sm_2(a,\dm I_d)-\frac{d}{2}.\ee
If the mask $a$ has a specific form, then we may be able to get a more accurate estimation of $\sm_\infty(a,\dm I_d)$. See Section~\ref{sec:exmp} for estimations of $\sm_\infty (a,4I_2)$ for specific two-dimensional masks $a$ in our examples.

\subsection{Motivation from CAGD and difficulties with stationary subdivision schemes}

In applications, we often require the underlying mask $a$ to satisfy certain critical properties for different purposes. First, one prefers a subdivision scheme that employs a mask with short support for the efficiency of implementation and computation. In the settings of CAGD, this is described by the size of \emph{stencils} of a mask (or a scheme). To be specific, for a mask $a\in\dlp{0}$, we define its \emph{filter support} to be the smallest $d$-dimensional interval $I:=[k_{1,1},k_{1,2}]\times [k_{2,1},k_{2,2}]\times \dots\times [k_{d,1},k_{d,2}]$ where $k_{1,1},k_{1,2},\dots,k_{d,1},k_{d,2}\in\Z$ such that $a(k)=0$ whenever $k\notin I$. A stationary $\dm I_d$-subdivision scheme with a mask $a\in\dlp{0}$ has \emph{$n$-ring stencils} for some $n\in\N$ if $\fsupp(a)\subseteq[-\dm n, \dm n]^d$. In general, it is highly desirable to have a $n$-ring stencil subdivision scheme such that $n$ is as small as possible, and that is, the mask $a$ has very small support. Next, a curve or a surface is modeled by a mesh by connecting neighborhood points. For example, in dimension two, there are two standard meshes: the triangular mesh and the quadrilateral mesh. A particular mesh is often associated with a symmetry group, and thus, the underlying mask of a subdivision scheme is often required to have the corresponding symmetry type.

In applications such as CAGD, people are particularly interested in a subdivision scheme that: (1) is interpolatory so that the limit function interpolates the initial data; (2) is at least $C^2$-convergent for the continuity of curvatures; (3) has no more than $2$-ring stencils to avoid exponentially increasing number of special subdivision rules near extraordinary vertices; (4) employs a mask with symmetry for a particular mesh. Unfortunately, these good properties cannot coexist in many cases. Suppose $a\in\dlp{0}$ is $\dm I_d$-interpolatory ($\dm\in\N\setminus\{1\}$), supported on $[-\dm,\dm]^d$ ( that is, has $1$-ring stencil) and satisfies $a(k)=a(-k)$ for all $k\in\dZ$ (this is the weakest symmetry type). On one hand, by \cite[Theorem 3.4]{hj06} (also see \cite[Theorem 4.1]{han00}), we have $\sr(a,2I_d)\le 2$ and therefore $\sm_\infty(a,\dm I_2)\le 2$. It then follows from \cite[Theorem 3.8]{hanphd} that the standard $\dm I_d$-refinable function $\phi$ of the mask $a$ is not a $C^2(\dR)$ function, so the stationary subdivision scheme that uses the mask $a$ cannot be $C^2$-convergent. Hence, any $C^2$-convergent interpolatory stationary $\dm I_d$-subdivision scheme that uses a symmetric mask must have at least $2$-ring stencils. On the other hand, for the most classical dilation matrix $2I_d$, it is pointed out in \cite[Corollary~4.3 and Theorem~3.5]{han00} and \cite[Theorem~3.9 and Corollary~3.12]{hanphd} that a $C^2$-convergent interpolatory stationary $2I_d$-subdivision scheme must have at least $3$-ring stencils. Consequently, the mutual conflict between properties (1)-(4) is a major difficulty with stationary subdivision schemes. Indeed, many existing famous stationary subdivision schemes do not satisfy all these properties. For instance, the famous butterfly scheme in \cite{dlg90} (also see \cite{cl23} for a non-linear analog of the scheme) and the interpolatory $2I_d$-subdivision schemes in \cite{hj98,rs97} all have no more than $2$-ring stencils but are not $C^2$-convergent; other modified butterfly schemes achieve $C^2$-convergence by either enlarging the support of the mask (\cite{nry16}) or making the mask to have $2$-ring stencils but sacrificing the interpolatory property (\cite{jy24}). Therefore, we need new settings and ideas to circumvent this difficulty with stationary subdivision schemes.

\subsection{Our contribution and paper structure}

To resolve the potential conflict between high-order smoothness and the short support of a refinement mask, following \cite{han23}, we introduce the notion of \emph{a quasi-stationary subdivision scheme}. Unlike a stationary subdivision scheme, a quasi-stationary subdivision scheme employs several different refinement masks repeatedly at different levels.
Let $\dm\in\N\setminus\{1\}$ be a dilation factor. For $r\in\N$ and $a_1,\dots,a_r\in\dlp{0}$, define
\[
\sd_{a_1,\dots,a_r,\dm I_d}^{n,r}:=\begin{cases}[\sd_{a_r,\dm I_d}\dots\sd_{a_1,\dm I_d}]^{\lfloor\frac{n}{r}\rfloor},&\text{if }n\in r\N_0,\\
\sd_{a_l,\dm I_d}\dots\sd_{a_1,\dm I_d}[\sd_{a_r,\dm I_d}\dots\sd_{a_1,\dm}]^{\lfloor\frac{n}{r}\rfloor},&\text{if }n\in r\N_0+l\text{ for some }l\in\{1,\dots,r-1\}.\end{cases}
\]
Suppose $\wh{a_l}(0)=1$ for all $l=1,\dots,r$. Given an initial data $v\in\dlp{0}$, the \emph{$r$-mask quasi-stationary $\dm I_d$-subdivision scheme} that uses $a_1,\dots,a_r$ generates a sequence $\{\sd_{a_1,\dots,a_r,\dm I_d}^{n,r}v\}_{n=1}^\infty$. When $r=1$, an $r$-mask quasi-stationary subdivision scheme becomes a stationary one. We define the $C^m$-convergence of a quasi-stationary subdivision scheme as the following:

\begin{definition}Let $\dm\in\N\setminus\{1\}$ be a dilation factor. Let $m\in\N_0$ and $a_1,\dots,a_r\in\dlp{0}$ be finitely supported filters such that $\wh{a_l}(0)=1$ for all $l=1,\dots,r$. If for every $v\in\dlp{0}$, there exists $\eta_v\in C^m(\dR)$ such that
	 \be\label{qss:def}\lim_{n\to\infty}\sup_{k\in\dZ}\left|\dm^{jn}[\nabla^\mu\sd_{a_1,\dots,a_r,\dm I_d}^{n,r}v](k)-\partial^\mu\eta_v(\dm^{-n}k)\right|=0,\quad\forall\mu\in\N_{0,j}^d ,\,j=0,\dots,m,\ee
	then we say that the $r$-mask quasi-stationary $\dm I_d$-subdivision scheme that uses $a_1,\dots,a_r$ is $C^m$-convergent.
	
\end{definition}

The main result of this paper is the following theorem, which fully characterizes the convergence and smoothness of a quasi-stationary subdivision scheme using only properties of the underlying masks.

\begin{theorem}\label{thm:qss}Let $\dm\in\N\setminus\{1\}$ be a dilation factor. Let $m\in\N_0$ and $a_1,\dots,a_r\in\dlp{0}$ be finitely supported filters such that $\wh{a_l}(0)=1$ for all $l=1,\dots,r$. Define $a\in\dlp{0}$ via
	 \be\label{qss:a}a:=\dm^{-dr}\sd_{a_r,\dm I_d}\dots\sd_{a_1,\dm I_d}\td.\ee
	The following statements are equivalent to each other:
	
	\begin{enumerate}
		
		\item[(1)]The $r$-mask quasi-stationary $\dm I_d$-subdivision scheme using $a_1,\dots,a_r$ is $C^m$-convergent;

		 \item[(2)]$\sm_\infty(a,\dm^rI_d)>m$ and $\sr(a_l,\dm I_d)>m$ for all $l=1,\dots,r$.
		
	\end{enumerate}
	If (1) or (2) holds, then for every $v\in\dlp{0}$, the limit function $\eta_v$ in \er{qss:def} must be given by
	 \be\label{eta:v}\eta_v=\sum_{k\in\Z^d}v(k)\phi(\cdot-k),\ee
	where $\phi$ is the standard $\dm^r I_d$-refinable function of the mask $a$:
	 \be\label{phi}\wh{\phi}(\xi):=\prod_{j=1}^{\infty}\wh{a}(\dm^{-rj}\xi),\quad\forall\xi\in\dR.\ee
	In particular, $\eta_{\td}=\phi$. Furthermore, the $r$-mask quasi-stationary $\dm I_d$-subdivision scheme that uses $a_1,\dots,a_r$ is interpolatory, that is,
\[
v(k)=\eta_v(k),\qquad\forall v\in\dlp{0},\quad k\in\dZ,
\]
	if and only if $a$ is $\dm^r I_d$-interpolatory, that is,
\[
a(\dm^rk)=\dm^{-rd}\td(k),\qquad\forall k\in\dZ.
\]
\end{theorem}

Let us comment about our contributions and explain the technicalities involved in our main result.

\begin{enumerate}
	\item[(1)] The special case $d=1$ of Theorem~\ref{thm:qss} has been established in \cite[Theorem~2 and Corollary~8]{han23}. As pointed out in \cite{han23}, the sum rule orders of the masks $a_1,\dots,a_r$ play a key role in analyzing the convergence and smoothness of a quasi-stationary subdivision scheme. In the case $d=1$, if a mask $u\in\lp{0}$ has order $m$ sum rules with respect to $\dm$, then $\wh{u}$ admits the following factorization:
\be\label{sr:d1}\wh{u}(\xi)=\left(1+e^{-i\xi}+\dots+e^{-(\dm-1)\xi}\right)^m\wh{b}(\xi),\quad\forall\xi\in\R,
\ee
for some finitely supported filter $b\in\lp{0}$. The factorization \er{sr:d1} is the key ingredient that greatly reduces the difficulty of the theoretical analysis in the case $d=1$. Unfortunately, a factorization like \er{sr:d1} is unavailable and often impossible when $d\ge 2$. Because of this, many tools from the case $d=1$ cannot be borrowed or directly generalized. We need new ideas to handle the case when $d$ is arbitrary.

	\item[(2)] In practice, one must estimate the $L_\infty$-smoothness exponent $\sm_\infty(a,\dm^rI_d)$ to analyze the smoothness order of a quasi-stationary subdivision scheme. In the case $d=1$, due to the simple characterization of the sum rule property in \er{sr:d1} of a one-dimensional filter, there are several simple and efficient methods to find lower estimates of $\sm_\infty(a,\dm^rI_d)$ (see \cite[Section 2.1]{han23} for a detailed survey). Unfortunately, tools from the one-dimensional case cannot be generalized to the multi-dimensional case in a straightforward way. Therefore, finding good estimates of $\sm_\infty(a,\dm^rI_d)$ is much more technical and difficult when $d\ge 2$.
As a consequence, except for tensor products of one-dimensional subdivision schemes, there are much fewer known multivariate subdivision schemes with high smoothness and small supports. See Section~\ref{sec:exmp} for detailed discussions on estimating the $L_\infty$-smoothness exponents of masks in our examples.

	\item[(3)] The dilation matrix $2I_d$ is the most classical choice and is most interesting in many applications. For $d=2$ and $r=2$, using the dilation matrix $2I_2$ of both triangular and quadrilateral meshes, for each smoothness exponent $m=1,2$, we provide in this paper examples of bivariate $C^m$-convergent interpolatory $2$-mask quasi-stationary $2I_2$-subdivision schemes such that all underlying masks $a_1,a_2$ have symmetry and at most $m$-ring stencils, i.e., $C^1$ smoothness with $1$-ring stencils and $C^2$-smoothness with $2$-ring stencils. These examples show that we can circumvent the difficulties with stationary subdivision schemes and demonstrate the advantages of quasi-stationary subdivision schemes.

\end{enumerate}

The structure of the paper is organized as follows:  In Section~\ref{sec:qss}, we first develop some auxiliary results regarding the sum rule properties in multi-dimensions. Then, we prove the main result Theorem~\ref{thm:qss}. In Section~\ref{sec:exmp}, we first briefly discuss constructing masks that satisfy all requirements of Theorem~\ref{thm:qss}. Next, we provide several illustrative examples of smooth interpolatory quasi-stationary subdivisions that use masks with, at most, $2$-ring stencils. We shall perform a detailed analysis on the $L_\infty$-smoothness exponent of the masks in our examples to prove the desired smoothness order of our schemes.

\begin{comment}
\textcolor{red}{the introduction is too short and it is a good idea to point out the contributions of this paper and the difficulties of the multivariate case.}
\end{comment}

\begin{comment}
\textcolor{red}{Using $\dm$ as a dilation factor is confusing for several places, e.g., you call $\phi$ is a $\dm$-refinable function. But $\dm$ is a scalar number and you are indicating that you are dealing with univariate case.
I think it is better to say something like: Throughout this paper, we consider a special $d\times d$ dilation matrix $\dm$ such that $dm=mI_d$ for some positive integer $m\ge 2$. But you also use $m$ for other purposes and this require you to change $m$ at other places. How about change the smoothness from $m$ to $s$ or $k$ or $n$ or $p$ etc. Another solution is to use $\dm I_d$ instead of $\dm$ to clearly indicate that you are dealing with dilation matrix $\dm I_d$ not dilation scalar factor $\dm$. For example, $\dm I_d$-subdivision scheme instead of $\dm$-subdivision}
\end{comment}

\section{Convergence and Smoothness of Quasi-stationary Subdivision Schemes}\label{sec:qss}

In this section, we prove the main result Theorem~\ref{thm:qss}.

\subsection{Auxiliary results}

To prove Theorem~\ref{thm:qss}, we must explore the sum rule properties of masks in $\dlp{0}$. To do this, we need the notion of \emph{coset masks}. Let $\dn$ be an invertible $d\times d$ integer matrix and define $d_{\dn}:=|\det(\dn)|$. For a mask $u\in\dlp{0}$ and $\gamma\in\dZ$, define the \emph{$\gamma$-coset mask of $u$ with respect to $\dn$} via
$$u^{[\gamma;\dn]}(k):=u(\gamma+\dn k),\qquad\forall k\in\dZ.$$
Using the definition of the Fourier series of $u$, it is easy to see that
\be\label{coset:u}\wh{u}(\xi):=\sum_{\gamma\in\Gamma_{\dn}}\wh{u^{[\gamma;\dn]}}(\dn^{\tp}\xi)e^{-i\gamma\cdot\xi},\qquad\forall \xi\in\dR,\ee
where $\Gamma_{\dn}$ is a complete set of representatives of the quotient group $\dZ/[\dn\dZ]$ and is given by
\be\label{ga:dn}\Gamma_{\dn}:=[\dn[0,1)^d]\cap\dZ:=\{\gamma_1,\dots,\gamma_{d_{\dn}}\}\text{ with }\gamma_1:=0.\ee
Define $\Omega_{\dn}$ to be a complete set of representatives of the quotient group $[\dn^{-\tp}\dZ]/\dZ$ given by
\be\label{om:dn}\Omega_{\dn}:=(\dn^{-\tp}\Z^d)\cap[0,1)^d:=\{\omega_1,\dots,\omega_{d_{\dn}}\}\text{ with }\omega_1:=0.\ee
It follows from \er{coset:u} that
\be\label{coset:0}[\wh{u}(\xi+2\pi\omega_1),\dots,\wh{u}(\xi+2\pi\omega_{d_{\dn}})]=[\wh{u^{[\gamma_{1};\dn]}}(\dn^\tp\xi),\dots,\wh{u^{[\gamma_{d_{\dn}};\dn]}}(\dn^\tp\xi)]F(\xi),\quad\forall\xi\in\dR,\ee
where $F(\xi)$ is the following $d_{\dn}\times d_{\dn}$ matrix:
\be\label{Fourier}F(\xi):=[e^{-i\gamma_j\cdot(\xi+2\pi\omega_l)}]_{1\le j,l\le d_{\dn}}.\ee
Noting that $F(\xi)\ol{F(\xi)}^\tp=d_{\dn}I_{d_{\dn}}$ for all $\xi\in\dR$, \er{coset:0} yields
\be\label{coset:1}\wh{u^{[\gamma_j;\dn]}}(\dn^\tp\xi)=d_{\dn}^{-1}e^{i\gamma_j\cdot\xi}\sum_{l=1}^{d_{\dn}}\wh{u}(\xi+2\pi\omega_l)e^{i\gamma_l\cdot(2\pi\omega_j)},\quad\forall j=1,\dots,d_{\dn},\quad \xi\in\dR.\ee

We have the following lemma.

\begin{lemma}\label{lem:sd1}Let $\dn$ be an invertible $d\times d$  integer matrix and define $\Omega_{\dn}$ via \er{om:dn}. Let $m\in\N_0$ and $u\in\dlp{0}$ be such that $\wh{u}(\xi+2\pi\omega)=\bo(\|\xi\|^m)$ as $\xi\to 0$ for all $\omega\in \Omega_{\dn}$, then
\be\label{sr:m}\wh{u}(\xi)=\sum_{\alpha\in\N_{0,m}^d}\wh{\nabla^\alpha\td}(\dn^\tp\xi)\wh{v_\alpha}(\xi),\qquad\forall\xi\in\dR,\ee
	for some $v_\alpha\in\dlp{0}$ for all $\alpha\in\N_{0,m}^d:=\{\nu\in\N_0^d:|\nu|=m\}$.
\end{lemma}

\bp Define $\Gamma_{\dn}$ as in \er{ga:dn}. By the assumption $\wh{u}(\xi+2\pi\omega)=\bo(\|\xi\|^m)$ as $\xi\to 0$ for all $\omega\in \Omega_{\dn}$ and \er{coset:1}, we have
$$\wh{u^{[\gamma_j;\dn]}}(\xi)=\bo(\|\xi\|^m),\quad\xi\to 0,\quad\forall j=1,\dots,d_{\dn}.$$
Hence, by \cite[Lemma 5]{dhacha}, we can write
$$\wh{u^{[\gamma_j;\dn]}}(\xi)=\sum_{\alpha\in\N_{0,m}^d}\wh{\nabla^\alpha\td}(\xi)\wh{u_{j,\alpha}}(\xi),\quad\forall\xi\in\dR,\quad j=1,\dots,d_{\dn},$$
for some mask $u_{j,\alpha}\in\dlp{0}$ for all $j\in\{1,\dots,d_{\dn}\}$ and $\alpha\in\N_{0,m}^d$. Therefore, we conclude from \er{coset:u} that \er{sr:m} must hold by choosing $v_\alpha\in\dlp{0}$ such that
$$\wh{v_\alpha}(\xi):=\sum_{j=1}^{d_{\dn}}\wh{u_{j,\alpha}}(\dn^\tp\xi)e^{-i\gamma_j\cdot\xi},\quad\forall\alpha\in\N_{0,m}^d,\quad\xi\in\dR.$$
This completes the proof.
\ep

\begin{remark}The special case of Lemma~\ref{lem:sd1} was proved in \cite[Lemma~2.5]{hanphd} with $\dn=2I_2$ and $m=1$ and the general case has already pointed out without proof in the remark after \cite[Theorem 3.6]{han03}.\end{remark}

With Lemma~\ref{lem:sd1}, we then have the following lemma on a crucial relation between the subdivision and the backward difference operators, which is essential to the proof of Theorem~\ref{thm:qss}.

\begin{lemma}\label{lem:sd2}Let $\dm\in\N\setminus\{1\}$, $m\in\N_0$ and $a\in\dlp{0}$ be such that $a$ has order $m$ sum rules with respect to $\dm I_d$. Then for every $\mu\in\N_{0,m}^d$, we have
	 \be\label{diff:sd}\nabla^\mu\sd_{a,\dm I_d}\td=\sum_{\alpha\in\N_{0,m}^d}\sd_{b_\alpha,\dm I_d}\nabla^\alpha\td,\ee
	 or equivalently,
	 \be\label{diff:sd:f}\wh{\nabla^\mu\td}(\xi)\wh{a}(\xi)=\sum_{\alpha\in\N_{0,m}^d}\wh{b_\alpha}(\xi)\wh{\nabla^\alpha\td}(\dm\xi),\quad \forall\xi\in\dR.\ee
	for some $b_\alpha\in\dlp{0}$ that satisfy
	 \be\label{sr:b}\wh{b_\alpha}(2\pi\omega)=\begin{cases}\td(\alpha-\mu){\dm}^{-m}\wh{a}(0),&\text{if }\omega=0,\\[0.2cm]
	 \frac{1}{(-i\dm)^m\alpha!}\wh{\nabla^\mu\td}(2\pi\omega)\partial^\alpha\wh{a}(2\pi\omega),&\text{otherwise},\end{cases}\quad \forall \alpha\in\N_{0,m}^d,\,\omega\in\Omega_{\dm}.
\ee
\end{lemma}

\bp All claims hold trivially if $m=0$, so we consider the case $m\in\N$. Since $a$ has order $m$ sum rules with respect to $\dm I_d$ and $\wh{\nabla^\mu\td}(\xi)=\bo(\|\xi\|^m)$ as $\xi\to 0$ for all $\mu\in\N_{0,m}^d$, it is clear that $\wh{\nabla^\mu\td}(\xi+2\pi\omega)\wh{a}(\xi+2\pi\omega)=\bo(\|\xi\|^m)$ as $\xi\to 0$ for all $\omega\in\Omega_{\dm}$. Hence, by Lemma~\ref{lem:sd1},  \er{diff:sd:f} must hold for some $b_\alpha\in\dlp{0}$ for all $\alpha\in\N_{0,m}^d$.

To prove \er{sr:b}, for every $\nu\in\N_{0,m}^d$, by taking the partial derivative $\partial^\nu$ on both sides of \er{diff:sd:f} and applying the product rule, we have
\be\label{diff:sd:1}\sum_{\beta\le\nu}\binom{\nu}{\beta}\partial^\beta\wh{\nabla^\mu\td}(\xi)\partial^{\nu-\beta}\wh{a}(\xi)=\sum_{\beta\le\nu}\sum_{\alpha\in\N_{0,m}^d}\binom{\nu}{\beta}\dm^{|\nu-\beta|}\partial^\beta\wh{b_\alpha}(\xi)\partial^{\nu-\beta}[\wh{\nabla^\alpha\td}](\dm\xi),\quad\forall\xi\in\dR.\ee
Now plug $\xi=0$ into \er{diff:sd:1}. Observe that $\partial^\beta\wh{\nabla^\mu\td}(0)$ is non-zero only if $\beta=\mu$ with $\partial^\mu\wh{\nabla^\mu\td}(0)=(-i)^m\mu!$, and $ \partial^{\nu-\beta}[\wh{\nabla^\alpha\td}](0)$ is non-zero only if $\nu=\alpha$ and $\beta=0$ with $\partial^{\alpha}[\wh{\nabla^\alpha\td}](0)=(-i)^m\alpha!$. Hence, by letting $\xi=0$ in \er{diff:sd:1} yields
$$\wh{b_\alpha}(0)=\td(\alpha-\mu)\dm^{-m}\wh{a}(0),\quad\forall \alpha\in\N_{0,m}^d.$$
Next, let $\omega\in\Omega_{\dm}\setminus\{0\}$ and plug $\xi=2\pi\omega$ into \er{diff:sd:1}. On one hand, since $a$ has order $m$ sum rules with respect to $\dm I_d$ and $\nu\in\N_{0,m}^d$, we have $\partial^{\nu-\beta}\wh{a}(2\pi\omega)=0$ for all $\beta\le \nu$ with $\beta\ne 0$. On the other hand, note that
$\partial^{\nu-\beta}[\wh{\nabla^\alpha\td}](\dm(2\pi\omega))$ is non-zero only if $\nu=\alpha$ and $\beta=0$ with $\partial^{\alpha}[\wh{\nabla^\alpha\td}](\dm(2\pi\omega))=(-i)^m\alpha!$. Hence, by letting $\xi=2\pi\omega$ with $\omega\in\Omega_{\dm}\setminus\{0\}$ in \er{diff:sd:1} yields
\[
\wh{b_\alpha}(2\pi\omega)=\frac{1}{(-i\dm)^m\alpha!}\wh{\nabla^\mu\td}(2\pi\omega)\partial^\alpha\wh{a}(2\pi\omega)\,\quad\forall\alpha \in\N_{0,m}^d,\,\omega\in\Omega_{\dm}\setminus\{0\}.
\]
The proof is now complete.
\ep

\begin{remark} The sum rule properties of multivariate masks have been investigated in the literature using other different approaches. For instance, Lemma~\ref{lem:sd1} was also investigated in \cite{sau02} for the special case $\dn=2I_d$ and then together with and the relation \er{diff:sd} in \cite{ms04,su08} for a general expansion matrix $\dn$ (i.e., all eigenvalues of $\dn$ are greater than $1$ in modulus). In the papers \cite{ms04,su08,sau02}, the authors characterize the sum rule properties from an algebraic perspective by using the theory of quotient ideals of Laurent polynomial rings, which requires a lot of prerequisites from algebra. Another possible approach to studying the sum rule properties is using the polynomial reproduction properties of the subdivision operator, see \cite{cc13,cry16} and many references therein. Our proof above follows the classical Fourier analytic method, which only uses the properties of Fourier series and coset masks. The Fourier analytic techniques give a more direct alternative approach to help us understand the sum rule properties and greatly facilitate the study of subdivision schemes.

\end{remark}

\begin{comment}
\textcolor{red}{You should mention that this results are also given in Sauer's paper. But here we present it as a special case of Lemma 2 and we explicitly construct $b_\alpha$ here.}

\end{comment}

\subsection{Proof of Theorem~\ref{thm:qss}}

We are now ready to prove Theorem~\ref{thm:qss}. We will first prove the more straightforward implication (2) $\Rightarrow$ (1) and then handle the more difficult implication (1) $\Rightarrow$ (2).

\bp[\textbf{Proof of Theorem~\ref{thm:qss}}] (2) $\Rightarrow$ (1):   Using linearity of the subdivision operator and the definition of the mask $a$ in \er{qss:a}, the quasi-stationary $\dm I_d$-subdivision operator that uses the masks $a_1,\dots,a_r$ is $C^m$ convergent if and only if
\be\label{qss:def:0}\lim_{n\to\infty}\|\dm^{|\mu|(rn+l)}\nabla^\mu\sd_{a_l,\dm I_d}\dots\sd_{a_1,\dm I_d}\sd_{a,\dm^rI_d}^n\td-\partial^\mu\eta_{\td}(\dm ^{-(rn+l)}\cdot)\|_{\ell_\infty(\Z^d)}=0,\quad\forall \mu\in\bigcup_{t=0}^m\N_{0,t}^d,
\ee
for all $l=0,1,\dots, r$, where $\eta_{\td}$ is the limit function of the particular input data $v=\td$. As item (2) holds, in particular $\sm_\infty(a,\dm^r I_d)>m$, we conclude from
\cite[Theorem 4.3]{han03} or
\cite[Theorem 7.3.1]{hanbook} that  \er{qss:def:0} holds with $l=0$ and $l=r$. Moreover, we must have $\eta_{\td}=\phi$ where $\phi$ is the standard $\dm^rI_d$-refinable function of $a$ that is defined as \er{phi}.

Next, we prove that \er{qss:def} must hold for $l=1,\dots,r-1$. Define $A_l\in\dlp{0}$ via
\be\label{A:l}\wh{A_l}(\xi):=\wh{a_1}(\dm^{l-1}\xi)\dots\wh{a_2}(\dm\xi)\wh{a_l}(\xi),\quad \forall\xi\in\dR,\,l=1,\dots,r-1.\ee
Since $\wh{a_l}(0)=1$ and $\sr(a_l,\dm I_d)>m$ for all $l=1,\dots,r$, we see that $\wh{A_l}(0)=1$ and $\sr(A_l,\dm^l I_d)> m$. In particular, for every $l=1,\dots,r-1$, we have
\be\label{sr:Al}\wh{A_l}(\xi+2\pi\omega_l)=\bo(\|\xi\|^{m+1}),\quad\xi\to 0,\quad \forall \omega_l\in\Omega_{\dm^l}:=[\dm^{-l}\Z^d]\cap[0,1)^d,\,\omega_l\ne 0.\ee
Let $j\in\{0,1,\dots,m\}$. For every $\mu\in\N_{0,j}^d$ and $l=1,\dots,r-1$, by Lemma~\ref{lem:sd2}, we can write
$$\wh{\nabla^\mu\sd_{a_l,\dm I_d}\dots\sd_{a_1,\dm I_d}\td}(\xi)=\wh{\nabla^\mu\td}(\xi)\wh{A_l}(\xi)=\sum_{\alpha\in\N_{0,j}^d}\wh{\nabla^\alpha\td}(\dm^l\xi)\wh{b_{l,\alpha}}(\xi)=\sum_{\alpha\in\N_{0,j}^d}\wh{\sd_{b_{l,\alpha},\dm^lI_d}\nabla^\alpha\td}(\xi),\quad\xi\in\dR,$$
for some $b_{l,\alpha}\in\dlp{0}$ for all $\alpha\in\N_{0,j}^d$ such that
\be\label{sr:bl}\wh{b_{l,\alpha}}(2\pi\omega_l)=\begin{cases}\td(\alpha-\mu)\dm^{-jl}, & \omega_l=0,\\
	0,&\omega_l\ne 0,\end{cases}\quad\forall \omega_l\in\Omega_{\dm^l}.\ee
For every $k\in\dZ$ and $n\in\N_0$, define
\be\label{jn}\begin{aligned}
	 J_n(k):=&\dm^{j(rn+l)}[\nabla^\mu\sd_{a_l,\dm I_d}\dots\sd_{a_1,\dm I_d}\sd_{a,\dm^rI_d}^n\td](k)-\partial^\mu\phi(\dm^{-(rn+l)}k)\\
	 =&\dm^{j(rn+l)}\left[\sum_{\alpha\in\N_{0,j}^d}\sd_{b_{l,\alpha},\dm^lI_d}\nabla^\alpha\sd_{a,\dm^rI_d}^n\td\right](k)-\partial^\mu\phi(\dm^{-(rn+l)}k).
\end{aligned}\ee
Then $J_{n}(k)=F_n(k)+G_n(k)$ where
\begin{align*}F_n(k):=&\dm^{jl}\sum_{\alpha\in\N_{0,j}^d}\left[\sd_{b_{l,\alpha},\dm^lI_d}\left(\dm^{jrn}\nabla^\alpha\sd_{a,\dm^rI_d}^n\td-\partial^\alpha\phi(\dm^{-rn}\cdot)\right)\right](k),\end{align*}
\begin{align*}
G_n(k):=&\dm^{jl}\left[\sum_{\alpha\in\N_{0,j}^d}\sd_{b_{l,\alpha},\dm^lI_d}(\partial^\alpha\phi(\dm^{-rn}\cdot))\right](k)-\partial^\mu\phi(\dm^{-(rn+l)}k)
\end{align*}
On the one hand, we have
\begin{align*}
\sup_{k\in\dZ}|F_n(k)|\le &\dm^{jl}\sum_{\alpha\in\N_{0,j}^d}\|\sd_{b_{l,\alpha},\dm^lI_d}\td\|_{\ell_1(\dZ)}\|\dm^{jrn}\nabla^\alpha\sd_{a,\dm^rI_d}^n\td-\partial^\alpha\phi(\dm^{-rn}\cdot)\|_{\ell_\infty(\dZ)}\\
\le &\dm^{(j+d)l}\sum_{\alpha\in\N_{0,j}^d}\|b_{l,\alpha}\|_{\ell_1(\dZ)}\|\dm^{jrn}\nabla^\alpha\sd_{a,\dm^rI_d}^n\td-\partial^\alpha\phi(\dm^{-rn}\cdot)\|_{\ell_\infty(\dZ)}.
\end{align*}
Since $\lim_{n\to\infty}\|\dm^{jrn}\nabla^\alpha[\sd_{a,\dm^rI_d}]^n\td-\partial^\alpha\phi(\dm^{-rn}\cdot)\|_{\ell_\infty(\dZ)}=0$ for all $\alpha\in\N_{0,j}^d$, we have
$$\lim_{n\to\infty}\sup_{k\in\dZ}|F_n(k)|=0.$$
On the other hand, by \er{sr:bl}, we conclude that
$$\sum_{q\in\dZ}b_{l,\alpha}(k+\dm^lq)=\dm^{-dl}\sum_{q\in\dZ}b_{l,\alpha}(q)=\td(\alpha-\mu)\dm^{-(j+d)l},\quad\forall \alpha\in\N_{0,j}^d,\,k\in\dZ.$$
It follows that
\begin{align*}
G_n(k)=&\dm^{(j+d)l}\sum_{\alpha\in\N_{0,j}^d}\sum_{q\in\dZ}b_{l,\alpha}(k-\dm^lq)\left[\partial^\alpha\phi(\dm^{-rn}q)-\partial^\alpha\phi(\dm^{-(rn+l)}k)\right]\\
=&\dm^{(j+d)l}\sum_{\alpha\in\N_{0,j}^d}\sum_{q\in \left(\dm^{-l}k-\dm^{-l}[-N,N]^d\right)\cap\dZ}b_{l,\alpha}(k-\dm^lq)\left[\partial^\alpha\phi(\dm^{-rn}q)-\partial^\alpha\phi(\dm^{-(rn+l)}k)\right],
\end{align*}
where $N\in\N$ is chosen such that $\cup_{\alpha\in\N_{0,j}^d}\fsupp(b_{l,\alpha})\subseteq [-N,N]^d$. For $q\in \left(\dm^{-l}k-\dm^{-l}[-N,N]^d\right)\cap\dZ$, we have
$$\|\dm^{-rn}q-\dm^{-(rn+l)}k\|=\dm^{-rn}\|q-\dm^{-l}k\|\le \dm^{-rn-l}\sqrt{d}N.$$
Hence
$$\sup_{k\in\dZ}|G_n(k)|\le \dm^{(j+d)l}\sum_{\alpha\in\N_{0,j}^d}\|b_{l,\alpha}\|_{\ell_1(\dZ)}\sup_{\|x-y\|\le \dm^{-rn-l}\sqrt{d}N}|\partial^\alpha \phi(x)-\partial^\alpha \phi(y)|.$$
Note that $\partial^\alpha\phi$ is compactly supported and uniformly continuous on $\R^d$ for all $\alpha\in\N_{0,j}^d$, we have
$$\lim_{n\to\infty}\sup_{\|x-y\|\le \dm^{-rn-l}\sqrt{d}N}|\partial^\alpha \phi(x)-\partial^\alpha \phi(y)|=0,$$
and thus
$$\lim_{n\to\infty}\sup_{k\in\dZ}|G_n(k)|=0.$$
Consequently,
$$\lim_{n\to\infty}\sup_{k\in\dZ}|J_n(k)|=0,$$
and this proves that \er{qss:def:0} holds for all $l=0,\dots,r$.

(1) $\Rightarrow$ (2): Suppose item (1) holds, that is, \er{qss:def} holds. By the definition of the mask $a$ in \er{qss:a}, we must have
\be\label{qss:def:1}\lim_{n\to\infty}\sup_{k\in\dZ}\left|\dm^{|\mu|rn}[\nabla^\mu\sd_{a,\dm^rI_d}^n\td](k)-\partial^\mu\eta_{\td}(\dm^{-rn}k)\right|=0,\quad\forall \mu\in\bigcup_{q=0}^m\N_{0,q}^d,\ee
By \cite[Theorem~4.3]{han03} or \cite[Theorem~7.3.1]{hanbook}, we must have $\sm_\infty(a,\dm^rI_d)>m$ and the limit function $\eta_{\td}=\phi$ must be the standard $\dm^r I_d$-refinable function associated with $a$ that is defined by \er{phi}.

Next, we show that $\sr(a_l,\dm I_d)>m$ for all $l=1,\dots,r$. Assume otherwise, that is, $\sr(a_l,\dm)=j\le m$ for some $l\in\{1,\dots,r\}$. For every $\mu\in\N_{0,j}^d$, by Lemma~\ref{lem:sd2}, we can write
\be\label{sr:al}\nabla^\mu\sd_{a_l,\dm I_d}=\sum_{\alpha\in\N_{0,j}^d}\sd_{b_{\mu,\alpha},\dm I_d}\nabla^\alpha,\ee
for some $b_{\mu,\alpha}\in\dlp{0}$ for all $\alpha\in\N_{0,j}^d$. For every $k\in\Z^d$ and $n\in\N_0$, define
$$J_n(k):=\dm^{j(rn+l)}[\nabla^\mu\sd_{a_l,\dm I_d}\dots\sd_{a_1,\dm}\sd_{a,\dm^rI_d}^n\td](k)-\partial^\mu\phi(\dm^{-(rn+l)}k).$$
We have $J_{n}(k)=H_n(k)+K_n(k)$ where
\begin{align*}H_n(k):=&\dm^{j}\sum_{\alpha\in\N_{0,j}^d}\left[\sd_{b_{\mu,\alpha},\dm I_d}\left(\dm^{j(rn+l-1)}\nabla^\mu\sd_{a_{l-1},\dm I_d}\dots\sd_{a,\dm^rI_d}^n\td-\partial^\alpha\phi(\dm^{-j(rn+l-1)}\cdot)\right)\right](k),\end{align*}
\begin{align*}
K_n(k):=&\dm^{j}\left[\sum_{\alpha\in\N_{0,j}^d}\sd_{b_{\mu,\alpha},\dm I_d}(\partial^\alpha\phi(\dm^{-j(rn+l-1)}\cdot))\right](k)-\partial^\mu\phi(\dm^{-j(rn+l)}k)
\end{align*}
By item (1), we have $\lim_{n\to\infty}\sup_{k\in\dZ}|J_n(k)|=0$ and
{\small\begin{align*}&\sup_{k\in\dZ}|H_n(k)|=\|H_n\|_{\ell_\infty(\dZ)}\\
	\le& \dm^{j-1}\sum_{\alpha\in\N_{0,j}^d}\left\|\sd_{b_{\mu,\alpha},\dm I_d}\left(\dm^{(j-1)(rn+l-1)}\nabla^\mu\sd_{a_{l-1},\dm I_d}\dots\sd_{a_1,\dm I_d}\sd_{a,\dm^rI_d}^n\td-\partial^\alpha\phi(\dm^{-(j-1)(rn+l-1)}\cdot)\right)\right\|_{\ell_\infty(\dZ)}\\
	\le &\dm^j\sum_{\alpha\in\N_{0,j}^d}\|b_{\mu,\alpha}\|_{\ell_1(\dZ)}\left\|\dm^{(j-1)(rn+l-1)}\nabla^\mu\sd_{a_{l-1},\dm I_d}\dots\sd_{a_1,\dm I_d}\sd_{a,\dm^rI_d}^n\td-\partial^\alpha\phi(\dm^{-(j-1)(rn+l-1)}\cdot)\right\|_{\ell_\infty(\dZ)}\\
	\to &0,\quad \text{ as } n\to\infty.\end{align*}   }
Therefore, we must have
\be\label{sup:g}\lim_{n\to\infty}\sup_{k\in\dZ}|K_n(k)|=\lim_{n\to\infty}\sup_{k\in\dZ}|J_n(k)-H_n(k)|=0.
\ee
By the definition of $\sd_{b_{\mu,\alpha},\dm I_d}$, we have
\begin{align*}
K_n(k)&=\dm^{j+d}\sum_{\alpha\in\N_{0,j}^d}\sum_{q\in\Z^d}b_{\mu,\alpha}(k-\dm q)\left[\partial^\alpha\phi(\dm^{-(rn+l-1)}q)-\partial^\mu\phi(\dm^{-(rn+l)}k)\right]\\
&=\dm^{j+d}\sum_{\alpha\in\N_{0,j}^d}\sum_{q\in\Z^d}b_{\mu,\alpha}(k-\dm q)\left[\partial^\alpha\phi(\dm^{-(rn+l-1)}q)-\partial^\alpha\phi(\dm^{-(rn+l-1)}k)\right]\\
\qquad &+\sum_{\alpha\in\N_{0,j}^d}c_{\mu,\alpha}(k)\partial^\alpha\phi(\dm^{-(rn+l-1)}k),
\end{align*}
where
\be\label{c:mu:alp}c_{\mu,\alpha}(k)=\dm^{j+d}\sum_{q\in\Z^d}b_{\mu,\alpha}(k-\dm q)-\dm^{j}\td(\alpha-\mu),\quad\forall \alpha\in\N_{0,d},\quad k\in\Z^d.\ee
Using the same argument as in the proof of $\lim_{n\to\infty}\sup_{k\in\Z^d}|G_n(K)|=0$ in the implication (2) $\Rightarrow$ (1), we can show that
\[
\lim_{n\to\infty}\sup_{k\in\Z^d}\left|\dm^{j+d}\sum_{\alpha\in\N_{0,j}^d}\sum_{q\in\Z^d}b_{\mu,\alpha}(k-\dm q)\left[\partial^\alpha\phi(\dm^{-(rn+l-1)}q)-\partial^\mu\phi(\dm^{-(rn+l)}k)\right]\right|=0,
\]
which, together with \eqref{sup:g} and identity on $K_n(k)$ after \eqref{sup:g}, forces
\be\label{phi:dev:0}\lim_{n\to\infty}\sup_{k\in\Z^d}\left|\sum_{\alpha\in\N_{0,j}^d}c_{\mu,\alpha}(k)\partial^\alpha\phi(\dm^{-(rn+l-1)}k)\right|=0.
\ee
Noting that $c_{\mu,\alpha}(k)=c_{\mu,\alpha}(k+\dm q_0)$ for all $k,q_0\in\dZ$, we conclude from \er{phi:dev:0} that
\be\label{phi:dev}\lim_{n\to\infty}\sup_{k\in \tilde{k}+\dm\Z^d}\left|\sum_{\alpha\in\N_{0,j}^d}c_{\mu,\alpha}(\tilde{k})\partial^\alpha\phi(\dm^{-(rn+l-1)}k)\right|=0,\quad \forall \tilde{k}\in\dZ.\ee
For every $\tilde{k}\in\dZ$, define
$$\Phi_{\tilde{k}}(x):=\sum_{\alpha\in\N_{0,j}^d}c_{\mu,\alpha}(\tilde{k})\partial^\alpha\phi(x),\quad\forall x\in\dR.$$
Let $x\in\dR$. For every $\eps>0$, by the uniform continuity of $\partial^\alpha\phi$ for all $\alpha\in\N_{0,j}^d$, we can choose $n_\eps\in\N_0$ and $z_\eps\in\dZ$ such that $|\Phi_{\tilde{k}}(x)-\Phi_{\tilde{k}}(\dm^{-n_\eps}z_{\eps})|<\eps$. For every $n\in\N_0$ such that $rn+l-1-n_\eps> 0$, we have
\begin{align*}
|\Phi_{\tilde{k}}(x)|&\le |\Phi_{\tilde{k}}(x)-\Phi_{\tilde{k}}(\dm^{-n_\eps}z_\eps)|+|\Phi_{\tilde{k}}(\dm^{-(rn+l-1)}\dm^{rn+l-1-n_\eps}z_\eps)-\Phi_{\tilde{k}}(\dm^{-(rn+l-1)}(\dm^{rn+l-1-n_\eps}z_\eps+\tilde{k}))|\\
&\qquad+|\Phi_{\tilde{k}}(\dm^{-(rn+l-1)}(\dm^{rn+l-1-n_\eps}z_\eps+\tilde{k}))|\\
&\le \eps+\sup_{\|s-t\|\le \dm^{-(rn+l-1)}\|\tilde{k}\|}|\Phi_{\tilde{k}}(s)-\Phi_{\tilde{k}}(t)|+\sup_{k\in\tilde{k}+\dm\Z^d}|\Phi_{\tilde{k}}(\dm^{-(rn+l-1)}k)|\\
&\to \eps,\text{ as }n\to\infty.
\end{align*}
Since $\eps>0$ is arbitrary, we must have $\Phi_{\tilde{k}}(x)=0$ for all $x\in\dR$.
 Now, taking the Fourier transform yields
$$0=\wh{\Phi_{\tilde{k}}}(\xi)=\sum_{\alpha\in\N_{0,j}^d}c_{\mu,\alpha}({\tilde{k}})\wh{[\partial^\alpha\phi]}(\xi)=\wh{\phi}(\xi)(-i)^j\sum_{\alpha\in\N_{0,j}^d}c_{\mu,\alpha}({\tilde{k}})\xi^\alpha,\quad\forall\xi\in\dR.$$
Note that $\phi$ has compact support with $\wh{\phi}(0)=1$, so $\wh{\phi}$ is a smooth analytic function which is not identically zero. This means we must have
$$\sum_{\alpha\in\N_{0,j}^d}c_{\mu,\alpha}({\tilde{k}})\xi^\alpha,\quad\forall\xi\in\dR.$$
Therefore
\be\label{c:mu:al}c_{\mu,\alpha}({\tilde{k}})=0,\quad\forall \alpha\in\N_{0,j}^d,\quad {\tilde{k}}\in\dZ.
\ee
By the definition of $c_{\mu,\alpha}(k)$ in \er{c:mu:alp}, \er{c:mu:al} is equivalent to say that $\sr(b_{\mu,\alpha},\dm I_d)\ge 1$ for all $\alpha\in\N_{0,j}^d$. In particular, we have $\wh{b_{\mu,\alpha}}(2\pi\omega)=0$ for all $\alpha\in\N_{0,j}^d$ and $\omega\in\Omega_{\dm}\setminus\{0\}$. On the other hand, by Lemma~\ref{lem:sd2}, we have
$$\wh{b_{\mu,\alpha}}(2\pi\omega)=\frac{1}{(-i\dm)^j\alpha!}\wh{\nabla^\mu\td}(2\pi\omega)\partial^\alpha\wh{a_l}(2\pi\omega),\quad\forall \alpha\in\N_{0,j}^d,\quad \omega\in\Omega_{\dm}\setminus\{0\},$$
which forces
$$\partial^\alpha\wh{a_l}(2\pi\omega)=0,\quad\forall \alpha\in\N_{0,j}^d,\quad \omega\in\Omega_{\dm}\setminus\{0\},$$
and thus implies that $\sr(a_l,\dm I_d)\ge j+1$, which is a contradiction. Therefore, the assumption $\sr(a_l,\dm I_d)\le m$ for some $l\in\{1,\dots,r\}$ is false and we must have $\sr(a_l,\dm I_d)>m$ for all $l=1,\dots,r$.

Consequently, items (1) and (2) must be equivalent to each other. The rest of the claims are trivial.\ep

\section{Examples of $2$-mask Interpolatory Quasi-stationary $2I_2$-Subdivision Schemes}\label{sec:exmp}

\begin{comment}
\textcolor{red}{You should use the technique in \cite{hj06} to perform the smoothness analysis. In particular, you need to prove that the Example 1 with $1$-ring is $C^1$ and Example 2 with $2$-ring is $C^2$. The analysis is not that hard and should be performed. In other words, the paper \cite{hj06} handles the dilation matrix $3I_2$ and here we will handle the case $4I_2$, which is very similar to thesis for dilation $4I_2$. If you read the maple code, then you can perform the analysis as in \cite{hj06}.
It is also a good idea to emphasize that for the other symmetry, one can obtain $C^1$ with 1-ring and $C^2$ with 2-ring by using tensor product. However, this paper we consider the other $D_6$-symmetry and one cannot use tensor product for this type of symmetry!}
\end{comment}

The dilation matrix $2I_d$ is of the most interest in the literature of subdivision schemes and wavelet theory. In this section, we present some examples for the case $d=2$ of $2$-mask interpolatory quasi-stationary $2 I_2$-subdivision schemes that use two symmetric masks $a_1$ and $a_2$. Moreover, our quasi-stationary subdivision schemes are $C^m$-convergent with $m$-ring stencils with $m\in\{1,2\}$. As previously mentioned, no $C^2$-convergent interpolatory stationary $2I_2$-subdivision scheme with two-ring stencils exists. Our example demonstrates that quasi-stationary subdivision schemes can overcome this shortcoming.

\subsection{Construction guideline}

We first discuss how to construct the masks $a_1,\dots,a_r$ that meet all requirements of Theorem~\ref{thm:qss}. Let $\dm\in\N\setminus\{1\}$ be a dilation factor and let $m\in\N_0$. We take the following steps to construct masks $a_1,\dots,a_r\in\dlp{0}$ such that $\wh{a_l}(0)=1$ for all $l=1,\dots,l$ and the $r$-mask interpolatory quasi-stationary $\dm I_d$-subdivision scheme that uses these masks is $C^m$-convergent:

\begin{enumerate}
	\item[(S1)] For each $l\in\{1,\dots,r\}$, parametrize the mask $a_l$ by
	$$\wh{a_l}(\xi)=\sum_{k\in [-\dm n,\dm n]^d}a_l(k)e^{-ik\cdot\xi},\qquad\forall j=1,\dots,r,\quad  \xi\in\dR,$$
for some $n\in\N$ so that $a_l$ has $n$-ring stencils for all $j=1,\dots,d$. Solve the linear system
	$$\sum_{k\in  [-\dm n,\dm n]^d}a_l(k)=1,\quad\forall l=1,\dots,r,$$
	and update $a_1,\dots,a_r$ by substituting in the solutions of the above system.
	
	\item[(S2)] \textbf{Sum rule conditions for $a_1,\dots a_r$:} For each $l\in\{1,\dots,r\}$, choose $n_l\in\N$ such that $n_l>m$. Solve the following linear system:
	 $$\sum_{k\in( [- n, n]^d-\dm^{-1}\gamma)\cap\dZ }(\gamma+\dm k)^\mu a_l(\gamma+\dm k)=\dm^{-d}\sum_{k\in  [-\dm n,\dm n]^d}k^\mu a_l(k),\quad  \gamma\in\Gamma_{\dm},\,\mu\in\N_0^d\text{ with }|\mu|<n_l,$$
	where $\Gamma_{\dm}=[0,\dm-1]^d\cap\dZ$. Update $a_1,\dots,a_r$ by substituting in the solutions of the above system.
	
	\item[(S3)] \textbf{Interpolatory condition:} Define $a\in\dlp{0}$ as in \er{qss:a}, or equivalently
	 $$\wh{a}(\xi):=\wh{a_1}(\dm^{r-1}\xi)\wh{a_2}(\dm^{r-2}\xi)\dots\wh{a_r}(\xi),\qquad\xi\in\dR.$$
	Solve the following system of equations:
	$$a(\dm^r k)=\dm^{-dr}\td(k),\quad \forall k\in [(\dm+\dm^{-1}+\dots+\dm^{-(r-1)})[-n,n]^d]\cap\dZ.$$
	Update $a, a_1,\dots,a_r$ by substituting in the solutions of the above system.
	
	\item[(S4)] \textbf{Try to optimize the smoothness exponent:} Choose the values of free parameters that make $\sm_2(a,\dm^rI_d)$, select parameter values among the remaining free parameters such that $\sm_2(a,\dm^rI_d)$ is as large as possible. Ideally, try to achieve $\sm_2(a,\dm^rI_d)>m+\frac{d}{2}$ so that $\sm_\infty(a,\dm^rI_d)>m$. If not possible, then try to directly estimate $\sm_\infty(a,\dm^rI_d)$ by using the structural properties of the mask $a$.
	
\end{enumerate}

By adding extra linear constraints to the above construction procedure, the masks $a_1,\dots,a_r$ that are constructed can also have \emph{symmetry}. The symmetry properties of multidimensional filters/masks are defined using the notion of symmetry groups. Let $\mathcal{G}$ be a finite set of $d\times d$ integer matrices that form a group under matrix multiplication. Here are some typically used symmetry groups in wavelet analysis:
\begin{itemize}
	
	\item $\mathcal{G}=\{-I_d,I_d\}$, where $I_d$ is the $d\times d$ identity matrix;
	
	\item For $d=2$, two important symmetry groups are
	
	\begin{itemize}
		
		\item \textbf{Full axis symmetry group:}
	\be\label{d4}D_4:=\left\{\pm I_2,\,\pm\begin{bmatrix}1 & 0\\
		0 & -1\end{bmatrix},\,\pm\begin{bmatrix}0 & 1\\
		1 & 0\end{bmatrix},\,\pm\begin{bmatrix}0 & 1\\
		-1 &0\end{bmatrix} \right\}.\ee
	This is the symmetry group associated with the quadrilateral mesh in $\Z^2$.
	
	\item \textbf{Hexagon symmetry group:}
	\be\label{d6}D_6:=\left\{\pm I_2,\,\pm\begin{bmatrix}0 & 1\\
		1 & 0\end{bmatrix},\,\pm\begin{bmatrix}-1 & 1\\
		0 & 1\end{bmatrix},\,\pm\begin{bmatrix}1 & 0\\
		1 & -1\end{bmatrix},\,\pm\begin{bmatrix}0 & 1\\
		-1 & 1\end{bmatrix},\, \pm\begin{bmatrix}1 & -1\\
		1 &0\end{bmatrix}\right\}.\ee
	This is the symmetry group associated with the triangular mesh in $\Z^2$.
	\end{itemize}
\end{itemize}

A filter $a\in\dlp{0}$ is \emph{$\mathcal{G}$-symmetric about a point $h\in\dR$} if
\be\label{sym:a}a(E(k-h)+h)=a(k),\qquad\forall k\in\dZ\text{ and }E\in\mathcal{G}.\ee
Let $\phi$ be the standard $\dm I_d$-refinable function associated with $a$ that is defined as \er{ref:a}. It is well-known that \er{sym:a} holds if and only if $\phi$ is \emph{$\mathcal{G}$-symmetric about a point $h_{\phi}:=(\dm-1)^{-1}h$}, that is,
\be\label{sym:phi}\phi(E(x-h_{\phi})+h_{\phi})=\phi(x),\qquad\forall x\in\dR\text{ and }E\in\mathcal{G}.\ee
If we require that all masks $a_1,\dots,a_r$ in an $r$-mask quasi-stationary subdivision scheme to have symmetry, then we can add the following linear constraint to the construction procedure above:
$$a_l(E(k-h_l)+h_l)=a_l(k),\qquad\forall l=1,\dots,r,\quad k\in\dZ,\quad E\in\mathcal{G},$$
for some selected points $h_l\in\dR$, $l=1,\dots,r$ and a given symmetry group $\mathcal{G}$.

For the case $d=2$ and $\dm=2$, people are interested in masks that are $D_4$- or $D_6$-symmetric in designing subdivision schemes or constructing wavelets and framelets. We will present examples of $C^m$-convergent $2$-mask interpolatory quasi-stationary $2I_2$-subdivision schemes with $m$-ring stencils for $m=1,2$. Particularly, unlike other types of symmetry, for $m\in\{1,2\}$, a $D_6$-symmetric mask with $m$-ring stencils that yields a $C^m$-convergence scheme cannot be obtained by performing tensor product of one-dimensional symmetric masks with $m$-ring stencils. Therefore,  $D_6$-symmetric examples are of significant importance and interest.

\subsection{Examples of $C^1$-convergent schemes} Let $r=2$ and $\dm=2$. We present two examples of $C^1$-convergent $2$-mask interpolatory quasi-stationary $2 I_2$-subdivision schemes with $1$-ring stencil using two symmetric masks $a_1$ and $a_2$.

For $u\in\dlpp{0}$, suppose $\fsupp(u):=[k_1,k_2]\times[n_1,n_2]$ for some $k_1,k_2,n_1,n_2\in\Z$. We use the following way to present a finitely supported filter $u\in l_0(\Z^2)$: suppose $\fsupp(u)=[k_1,k_2]\times[n_1,n_2]$, then we write
$$u=\begin{bmatrix}u(k_1,n_2) & u(k_1+1,n_2) &  \dots & u(k_2,n_2)\\
u(k_1,n_2-1) & u(k_1+1,n_2-1) &  \dots & u(k_2,n_2-1)\\
\vdots & \vdots &\ddots &\vdots\\
u(k_1, n_1) & u(k_1+1,n_1) &  \dots & u(k_2,n_1)\end{bmatrix}_{[k_1,k_2]\times[n_1,n_2]}.$$
For example, $\wh{u}(\xi_1,\xi_2)=e^{-i\xi_1}+2e^{i\xi_2}$ is presented as $u=\begin{bmatrix}
0 & 1\\
2 &0
\end{bmatrix}_{[0,1]\times[-1,0]}$.

Let $\mathcal{G}=D_4$ and parameterize two masks $a_1,a_2\in l_0(\Z^2)$ such that \begin{itemize}
	\item $\wh{a_1}(0)=\wh{a_2}(0)=1$, $\sr(a_1,2I_2)=\sr(a_2,2I_2)=2$;
	
	\item $a_1$ and $a_2$ have $1$-ring stencil;
	
	\item $a_1$ and $a_2$ are $D_4$-symmetric about $(0,0)$;
\end{itemize}
as follows:
\be\label{1:a1:d4}\wh{a_1}(\xi_1,\xi_2):=\frac{1}{16}e^{i(\xi_1+\xi_2)}
(1+e^{-i\xi_1})^2(1+e^{-i\xi_2})^2\wh{q_1}(\xi_1,\xi_2),
\ee
\be\label{1:a2:d4}\wh{a_2}(\xi_1,\xi_2):=\frac{1}{16}e^{i(\xi_1+\xi_2)}
(1+e^{-i\xi_1})^2(1+e^{-i\xi_2})^2\wh{q_2}(\xi_1,\xi_2),
\ee
where $q_1,q_2\in l_0(\Z^2)$ are given by
\be\label{d4:q}q_1=\begin{bmatrix}t_2  & t_1 & t_2\\
t_1 & 1-4t_1-4t_2 & t_1\\
t_2  & t_1 & t_2\end{bmatrix}_{[-1,1]^2},\quad q_2=\begin{bmatrix}t_4  & t_3 & t_4\\
t_3 & 1-4t_3-4t_4 & t_3\\
t_4  & t_3 & t_4\end{bmatrix}_{[-1,1]^2}, \ee
for some free parameters $t_1,t_2,t_3,t_4\in\R$. Define the mask $a\in l_0(\Z^2)$  by
\be\label{1:a}a:=2^{-4}\sd_{a_2,2I_2}\sd_{a_1,2I_2}\td,\ee
or equivalently
$$\wh{a}(\xi_1,\xi_2):=\wh{a_1}(2\xi_1,2\xi_2)\wh{a_2}(\xi_1,\xi_2).$$ We must have $\sm_\infty(a,4I_2)>1$ to guarantee the $C^1$-convergence of the $2$-mask quasi-stationary subdivision scheme. To achieve this, we need to have a reasonable estimation of $\sm_\infty(a,4I_2)$ and then choose the values of the free parameters that yield the desired result. We recall the following results from \cite{hj06}:

\begin{theorem}\label{thm:sm:1}(\cite[Theorem 2.4]{hj06})Let $\dm\in\N\setminus\{1\}$ and let $b\in\dlp{0}$ be a filter. Then
	 \be\label{est:rho:1}\lim_{n\to\infty}\left\|\sd^n_{b,\dm I_d}\td]\right\|_{l_\infty(\dZ)}^{\frac{1}{n}}=\inf_{n\in\N}\left(\sup_{\gamma\in\Gamma_{\dm}}\sum_{k\in\dZ}\left|\sd_{b,\dm I_d}^n\td(\gamma+\dm^nk)\right|\right)^{\frac{1}{n}},\ee
	where $\Gamma_{\dm}:=[0,\dm-1]^d\cap\dZ$.
\end{theorem}

\begin{theorem}\label{thm:sm:2}
(\cite[Theorem 2.3]{hj06})
Let $\dm\in\N\setminus\{1\}$ and let $u\in\dlp{0}$ be a filter. If $\wh{u}(\xi)=\frac{\wh{c}(\dm\xi)}{\wh{c}(\xi)}\wh{b}(\xi)$ for some filters $b,c\in\dlp{0}$ such that $\frac{\wh{c}(\dm\xi)}{\wh{c}(\xi)}$ is a $2\pi\dZ$-periodic trigonometric polynomial, then
	 \be\label{est:rho}\lim_{n\to\infty}\left\|v*c*[\sd^n_{u,\dm I_d}\td]\right\|_{l_\infty(\dZ)}^{\frac{1}{n}}=\lim_{n\to\infty}\left\|u*[\sd^n_{b,\dm I_d}\td]\right\|_{l_\infty(\dZ)}^{\frac{1}{n}},\qquad \forall v\in\dlp{0}.\ee
\end{theorem}

Let $a\in l_0(\Z^2)$ be defined by \er{1:a}, where $a_1,a_2\in l_0(\Z^2)$ are defined by \er{1:a1:d4} and \er{1:a2:d4}. To estimate $\sm_\infty(a,4I_2)$, it suffices to estimate
$$\rho_2(a,4I_2,\mu)_\infty:=\lim_{n\to\infty}\left\|\nabla^\mu\sd_{a,4I_2}^n\td\right\|_{e_\infty(\Z^2)}^{\frac{1}{n}},\quad\text{for }\mu=(2,0),\,(1,1),\,(0,2).$$
Since $a_1,a_2$ are $D_4$-symmetric about $(0,0)$ and $\sr(a_1,2I_2)=\sr(a_2,2I_2)=2$, so is the mask $a$, and $\sd_{a,4I_2}^n\td$ is also $D_4$-symmetric about $(0,0)$. For any mask $u\in l_0(\Z^2)$ that is $D_4$-symmetric about $(0,0)$, we have
\be\label{d4:u}\nabla^{(2,0)} u(Ek)=\nabla_{e_1}^2u(Ek)=\nabla_{E^{-1}e_1}^2u(k),\quad\forall k\in\Z^2.\ee
By letting $E=\begin{bmatrix}0 &1\\
1 &0\end{bmatrix}\in D_4$, \er{d4:u} yields
\be\label{d4:u:2}\nabla^{(2,0)} u(Ek)=\nabla_{e_1}^2u(Ek)=\nabla^{(0,2)}u(k),\quad\forall k\in\Z^2.\ee
It then follows from \er{d4:u:2} that $\rho_2(a,4I_2,(2,0))_\infty=\rho_2(a,4I_2,(0,2))_\infty$ and thus
\be\label{rho:a:d4}\rho_2(a,4I_2)_\infty=\max\{\rho_2(a,4I_2,(2,0))_\infty,\,\rho_2(a,4I_2,(1,1))_\infty\}.
\ee
Let $q_1,q_2\in l_0(\Z^2)$ be defined as in \er{d4:q}. Define $b_1,b_2\in l_0(\Z^2)$ via
\be\label{d4:b1}\wh{b_1}(\xi_1,\xi_2)=\frac{1}{256}e^{3i(\xi_1+\xi_2)}
(1+e^{-i\xi_2})^2(1+e^{-2i\xi_2})^2\wh{q_1}(2\xi_1,2\xi_2)\wh{q_2}(\xi_1,\xi_2),\ee
\be\label{d4:b2}\wh{b_2}(\xi_1,\xi_2)=\frac{1}{256}e^{3i(\xi_1+\xi_2)}
(1+e^{-i\xi_1})(1+e^{-i\xi_2})(1+e^{-2i\xi_1})(1+e^{-2i\xi_2})\wh{q_1}(2\xi_1,2\xi_2)\wh{q_2}(\xi_1,\xi_2).
\ee
Note that
\[
\wh{a}(\xi_1,\xi_2)=\frac{\wh{\nabla^{(2,0)}\td}(4\xi_1,4\xi_2)}{\wh{\nabla^{(2,0)}\td}(\xi_1,\xi_2)}\wh{b_1}(\xi_1,\xi_2)=
\frac{\wh{\nabla^{(1,1)}\td}(4\xi_1,4\xi_2)}{\wh{\nabla^{(1,1)}\td}(\xi_1,\xi_2)}\wh{b_2}(\xi_1,\xi_2).
\]
It then follows from Theorems~\ref{thm:sm:1} and ~\ref{thm:sm:2} that
\begin{align*}
&\rho_2(a,4I_2,(2,0))_\infty=\lim_{n\to\infty}\left\|\sd_{b_1,4I_2}^n\td\right\|_{l_\infty(\Z^2)}^{\frac{1}{n}}\le\left(\sup_{\gamma\in\Gamma_{4^n}}\sum_{k\in\Z^2}\left|\sd_{b_1,4I_2}^n \td(\gamma+4^nk)\right|\right)^{\frac{1}{n}},\\
&\rho_2(a,4I_2,(1,1))_\infty=\lim_{n\to\infty}\left\|\sd_{b_2,4I_2}^n\td\right\|_{l_\infty(\Z^2)}^{\frac{1}{n}}\le\left(\sup_{\gamma\in\Gamma_{4^n}}\sum_{k\in\Z^2}\left|\sd_{b_2,4I_2}^n \td(\gamma+4^nk)\right|\right)^{\frac{1}{n}},
\end{align*}
for all $n\in\N$ where $\Gamma_{4^n}:=[0,4^n-1]^2\cap\Z^2$. Therefore, we have
\be\label{est:sm:1}\sm_\infty(a,4I_2)\ge-\log_4\left(\max_{j=1,2}\left(\sup_{\gamma\in\Gamma_{4^n}}\sum_{k\in\Z^2}\left|\sd_{b_j,4I_2}^n \td(\gamma+4^nk)\right|\right)^{\frac{1}{n}}\right),\quad\forall n\in\N,
\ee
where $b_1,b_2\in l_0(\Z^2)$ are given by \er{d4:b1} and \er{d4:b2}.

\begin{figure}[htbp]
	\centering
	\begin{subfigure}[b]{0.24\textwidth} \includegraphics[width=\textwidth,height=0.8\textwidth]{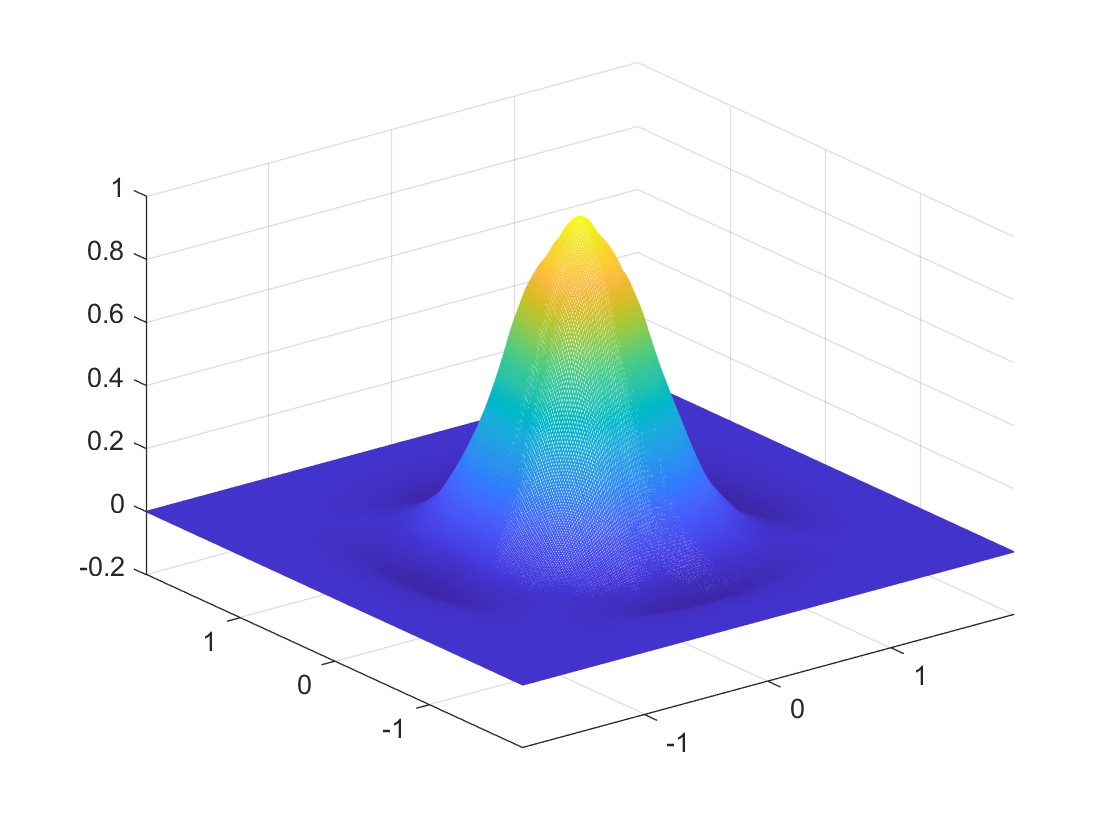}
		\caption{$\phi\in \mathscr{C}^1(\R^2)$}
	\end{subfigure}
	\begin{subfigure}[b]{0.22\textwidth} \includegraphics[width=\textwidth,height=0.8\textwidth]{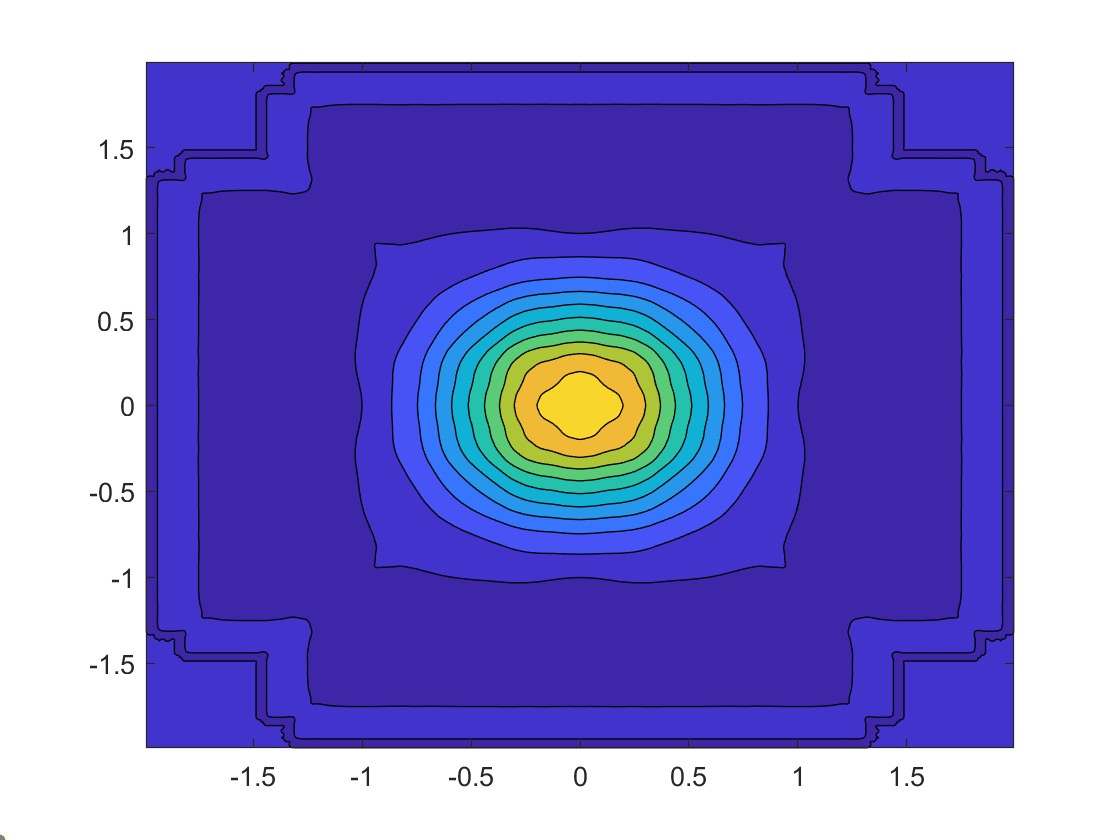}
		\caption{Contour of $\phi$}
	\end{subfigure}
	\begin{subfigure}[b]{0.24\textwidth} \includegraphics[width=\textwidth,height=0.8\textwidth]{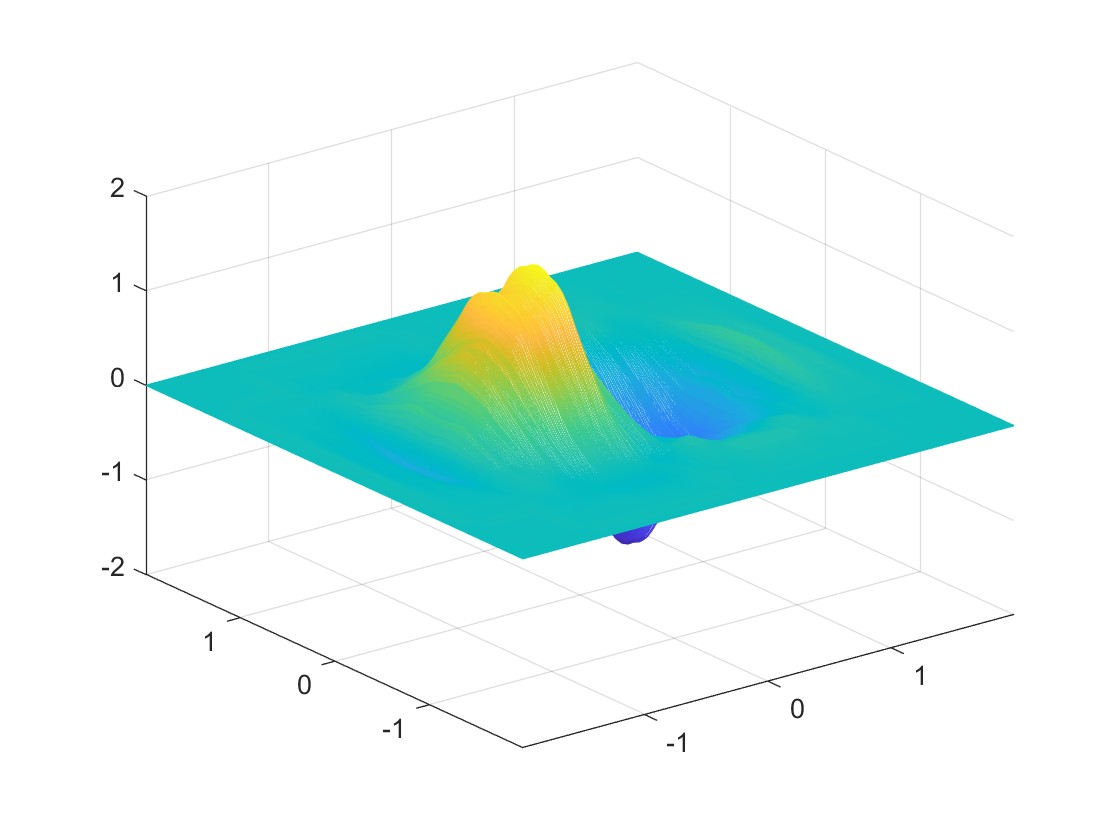} \caption{$\frac{\partial \phi}{\partial x} \in C^0(\R^2)$}
	\end{subfigure}
	\begin{subfigure}[b]{0.22\textwidth} \includegraphics[width=\textwidth,height=0.8\textwidth]{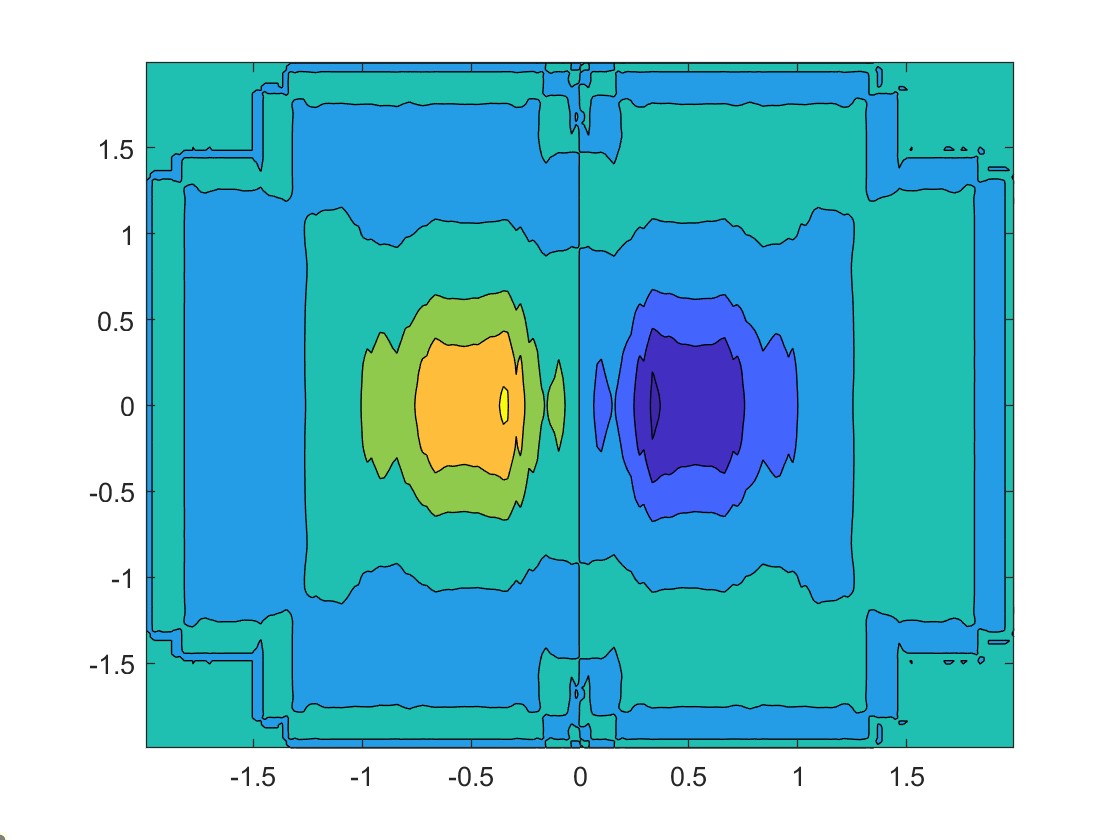}
		\caption{Contour of $\frac{\partial \phi}{\partial x}$}
	\end{subfigure}
	\caption{
		(A) is the graph of the interpolating $4I_2$-refinable function $\phi\in C^1(\R^2)$ in Example~\ref{ex3} and (B) is its contour.
		(C) is the graph of the partial derivative $\frac{\partial \phi}{\partial x} \in C^0(\R^2)$, and (D) is the contour of $\frac{\partial \phi}{\partial x}$.
	}\label{fig:ex3}
\end{figure}

\begin{exmp}\label{ex3} Let $q_1,q_2\in l_0(\Z^2)$ be given by \er{d4:q} for some free parameters $t_1,t_2,t_3,t_4\in\R$. Let $a_1,a_2\in l_0(\Z^2)$ be given by \er{1:a1:d4} and \er{1:a2:d4}. Define $a\in l_0(\Z^2)$ via \er{1:a}. By imposing $t_2=0$ and the $4I_2$-interpolatory constraint $a(4k)=\frac{1}{16}\td(k)$ for all $k\in\Z^2$, we obtain
\[
t_1=t_1,\quad t_2=0,\quad t_3=-\tfrac{2t_1(4t_1-1)}{8t_1^2-6t_1+1},\quad t_4=\tfrac{8t_1^2}{8t_1^2-6t_1+1}.
\]
By taking $t_1=-\frac{11}{42}$, in which case $t_3=-\frac{11}{1376}$ and $t_4=\frac{121}{688}$,  the two masks $a_1,a_2$ are then given by
$$a_1=\frac{1}{672}\begin{bmatrix}0 & -11 & -22 & -11 &  0\\
-11 & 42 & 106 & 42 & -11\\
-22 & 106 & 256 & 106 & -22\\
-11 & 42 & 106 & 42 & -11\\
0 &  -11 &  -22 & -11 &  0\end{bmatrix}_{[-2,2]^2},\quad a_2=\frac{1}{22016}\begin{bmatrix}242 & 473 & 462 & 473 & 242\\
473 & 1376 & 1806 & 1376 & 473\\
462 & 1806 & 2688 & 1806 & 462\\
473 & 1376 & 1806 & 1376 & 473\\
242 & 473 & 462 & 473 & 242\end{bmatrix}_{[-2,2]^2}.$$
Computation yields  $\sm_2(a,4I_2)\approx1.70906$. Using the estimation in \er{est:sm:1} with $n=2$, we obtain $\sm_\infty(a,4I_2)\ge1.38616$. Therefore, the $2$-mask interpolatory quasi-stationary $2I_2$-subdivision scheme that uses the above masks $a_1,a_2$ is $C^1$-convergent. See Figure~\ref{fig:ex3} for the graphs of the standard $4I_2$-refinable function $\phi$ of the mask $a$ , $\frac{\partial \phi}{\partial x}$ and the contours of $\phi$ and $\frac{\partial \phi}{\partial x}$.

\begin{comment}
Another choice is
\[
t_1=-\frac{88}{105},\quad t_2=\frac{121}{420}, \quad
t_2=\frac{11}{32},\quad t_4=0.
\]
\end{comment}

\end{exmp}

Let $\mathcal{G}=D_6$ and parameterize two masks $a_1,a_2\in l_0(\Z^2)$ such that \begin{itemize}
	\item $\wh{a_1}(0)=\wh{a_2}(0)=1$, $\sr(a_1,2I_2)=\sr(a_2,2I_2)=2$;
	
	\item $a_1$ and $a_2$ have $1$-ring stencil;
	
	\item $a_1$ and $a_2$ are $D_6$-symmetric about $(0,0)$;
\end{itemize}
as follows
\be\label{1:a1}\wh{a_1}(\xi_1,\xi_2):=\frac{1}{8}(1+e^{-i\xi_1})(1+e^{-i\xi_2})(1+e^{i(\xi_1+\xi_2)})(t_1\wh{p}(\xi_1,\xi_2)+1),
\ee
\be\label{1:a2}\wh{a_2}(\xi_1,\xi_2):=\frac{1}{8}(1+e^{-i\xi_1})(1+e^{-i\xi_2})(1+e^{i(\xi_1+\xi_2)})(t_2\wh{p}(\xi_1,\xi_2)+1),
\ee	
where $t_1,t_2$ are free parameters and $p\in l_0(\Z^2)$ is given by
\be\label{1:p}\wh{p}(\xi):=e^{-i\xi_1}+e^{i\xi_1}+e^{-i\xi_2}+e^{i\xi_2}+e^{-i(\xi_1+\xi_2)}+e^{i(\xi_1+\xi_2)}-6.\ee
Define the mask $a\in l_0(\Z^2)$ via \er{1:a} with the above masks $a_1,a_2\in l_0(\Z^2)$ and we estimate $\sm_\infty(a,4I_2)$. Since \er{d4:u} must hold for all $u\in l_0(\Z^2)$ and $E=\begin{bmatrix} 0 & 1\\ 1 & 0\end{bmatrix}\in D_6$, it is easy to see that \er{rho:a:d4} must hold.  For every $u\in l_0(\Z^2)$ and $k\in\Z^2$, we have
\be\label{n:e1}\nabla_{e_1}u(k)=u(k)-u(k+e_2)+u(k+e_2)-u(k-e_1)=-\nabla_{e_2}u(k+e_2)+\nabla_{e_1+e_2}u(k+e_2).\ee
Hence
\begin{align*}\left\|\nabla^{(2,0)}\sd_{a,4I_2}^n\td\right\|_{l_\infty(\Z^2)}&=\left\|\nabla_{e_1}^2\sd_{a,4I_2}^n\td\right\|_{l_\infty(\Z^2)}\le \left\|\nabla_{e_1}\nabla_{e_2}\sd_{a,4I_2}^n\td\right\|_{l_\infty(\Z^2)}+\left\|\nabla_{e_1}\nabla_{e_1+e_2}\sd_{a,4I_2}^n\td\right\|_{l_\infty(\Z^2)}\\
&\le 2\max\left\{\left\|\nabla_{e_1}\nabla_{e_2}\sd_{a,4I_2}^n\td\right\|_{l_\infty(\Z^2)},\,\left\|\nabla_{e_1}\nabla_{e_1+e_2}\sd_{a,4I_2}^n\td\right\|_{l_\infty(\Z^2)}\right\},
\end{align*}
and thus
$$\rho_2(a,4I_2,(2,0))_\infty\le\max\left\{\lim_{n\to\infty}\left\|\nabla_{e_1}\nabla_{e_2}\sd_{a,4I_2}^n\td\right\|_{l_\infty(\Z^2)}^{\frac{1}{n}},\,\lim_{n\to\infty}\left\|\nabla_{e_1}\nabla_{e_1+e_2}\sd_{a,4I_2}^n\td\right\|_{l_\infty(\Z^2)}^{\frac{1}{n}}\right\}.$$
Moreover, as $E=\begin{bmatrix}1 &-1\\ 0 & -1\end{bmatrix}\in D_6$, the $D_6$-symmetry of $\sd_{a,4I_2}^n\td$ yields
\be\label{e1:e2}\left\|\nabla_{e_1}\nabla_{e_2}\sd_{a,4I_2}^n\td\right\|_{l_\infty(\Z^2)}=\left\|\nabla_{E^{-1}e_1}\nabla_{E^{-1}e_2}\sd_{a,4I_2}^n\td\right\|_{l_\infty(\Z^2)}=\left\|\nabla_{e_1}\nabla_{e_1+e_2}\sd_{a,4I_2}^n\td\right\|_{l_\infty(\Z^2)}.\ee
Therefore, it follows that
$$\rho_2(a,4I_2,(2,0))_\infty\le \lim_{n\to\infty}\left\|\nabla_{e_1}\nabla_{e_2}\sd_{a,4I_2}^n\td\right\|_{l_\infty(\Z^2)}^{\frac{1}{n}}=\rho_2(a,4I_2,(1,1))_\infty,$$
which further implies that
$$\rho_2(a,4I_2)_\infty=\rho_2(a,4I_2,(1,1))_\infty.$$
Let $p\in l_0(\Z^2)$ be the same as in \er{1:p}. Define $b\in l_0(\Z^2)$ via
\be\label{b}\wh{b}(\xi):=\frac{1}{64}(1+e^{i(\xi_1+\xi_2)})(1+e^{2i(\xi_1+\xi_2)})(t_1\wh{p}(2\xi)+1)(t_2\wh{p}(\xi)+1).\ee
Note that
$$\wh{a}(\xi_1,\xi_2)=\frac{\wh{\nabla^{(1,1)}\td}(4\xi_1,4\xi_2)}{\wh{\nabla^{(1,1)}\td}(\xi_1,\xi_2)}\wh{b}(\xi_1,\xi_2),$$
We then conclude from Theorems~\ref{thm:sm:1} and ~\ref{thm:sm:2} that
\[
\rho_2(a,4I_2,(1,1))_\infty=\lim_{n\to\infty}\left\|\sd_{b,4I_2}^n\td\right\|_{l_\infty(\Z^2)}^{\frac{1}{n}}\le \left(\sup_{\gamma\in\Gamma_{4^n}}\sum_{k\in\dZ}\left|\sd_{b,4I_2}^n\td(\gamma+\dm^nk)\right|\right)^{\frac{1}{n}},\quad \forall n\in\N.
\]
Consequently, we obtain
\be\label{sm:a:1}
\sm_\infty(a,4I_2)\ge-\log_4\left(
\sup_{\gamma\in\Gamma_{4^n}}\sum_{k\in\Z^2}
\left|\sd_{b,4I_2}^n\td(\gamma+4^nk)\right|
\right)^{\frac{1}{n}},\qquad\forall n\in\N,
\ee
where $b\in l_0(\Z^2)$ is given by \er{b}.

\begin{figure}[htbp]
	\centering
	\begin{subfigure}[b]{0.24\textwidth} \includegraphics[width=\textwidth,height=0.8\textwidth]{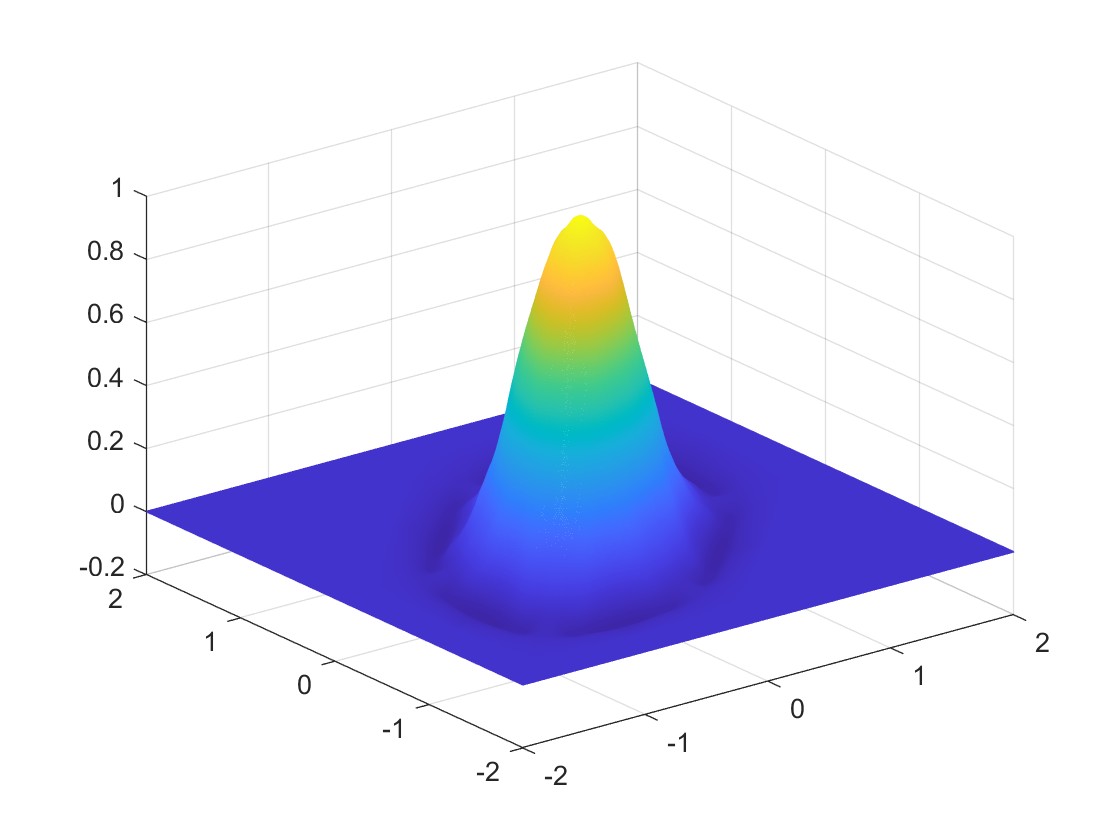}
		\caption{$\phi\in \mathscr{C}^1(\R^2)$}
	\end{subfigure}
	\begin{subfigure}[b]{0.22\textwidth} \includegraphics[width=\textwidth,height=0.8\textwidth]{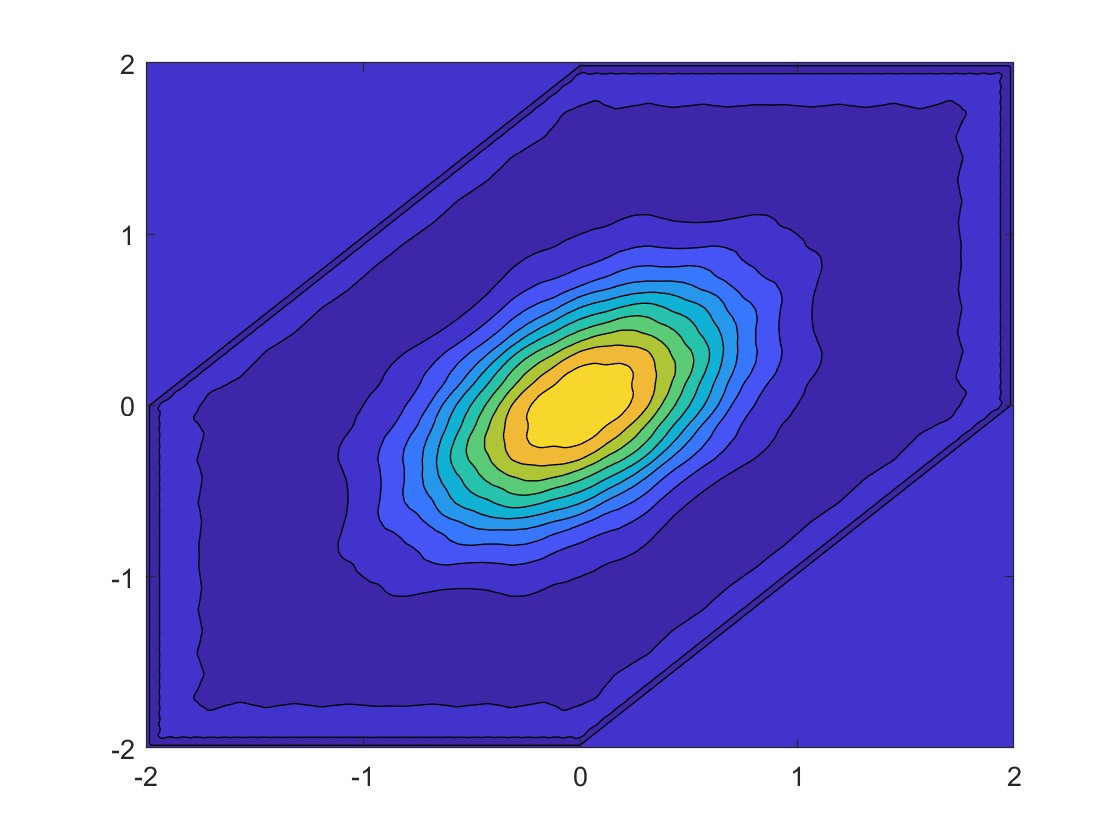}
		\caption{Contour of $\phi$}
	\end{subfigure}
	\begin{subfigure}[b]{0.24\textwidth} \includegraphics[width=\textwidth,height=0.8\textwidth]{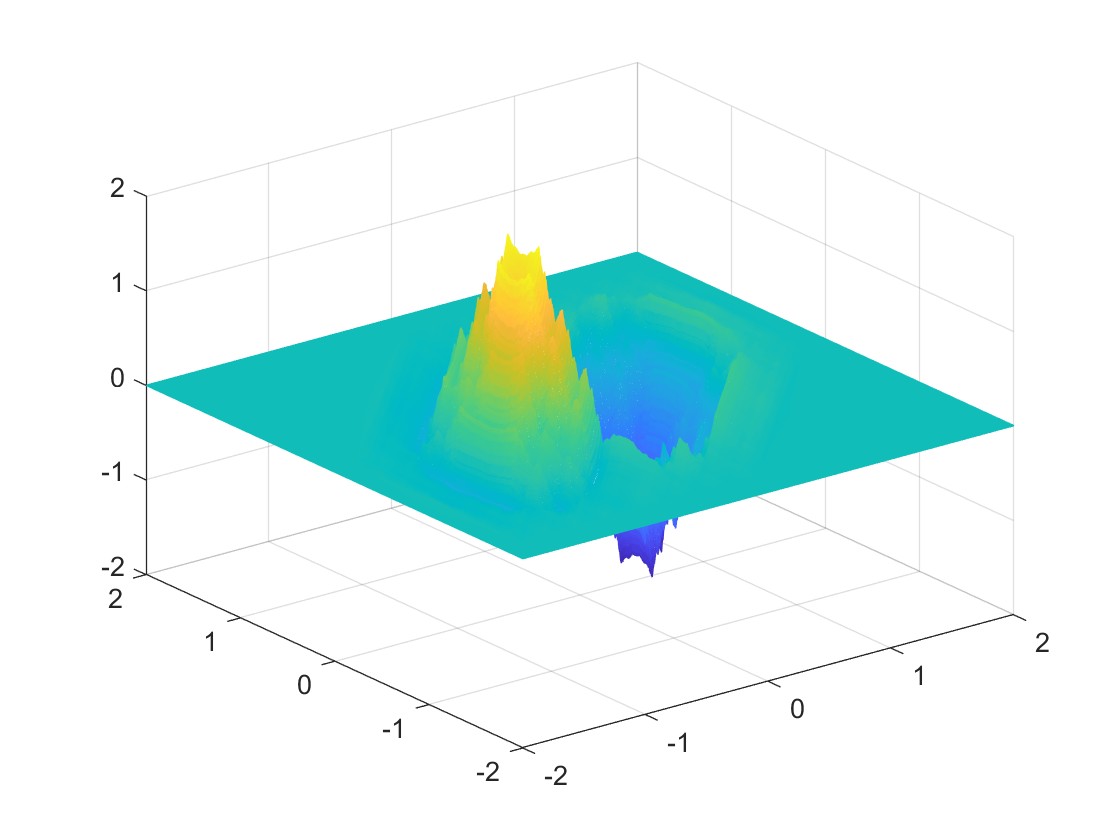} \caption{$\frac{\partial \phi}{\partial x} \in C^0(\R^2)$}
	\end{subfigure}
	\begin{subfigure}[b]{0.22\textwidth} \includegraphics[width=\textwidth,height=0.8\textwidth]{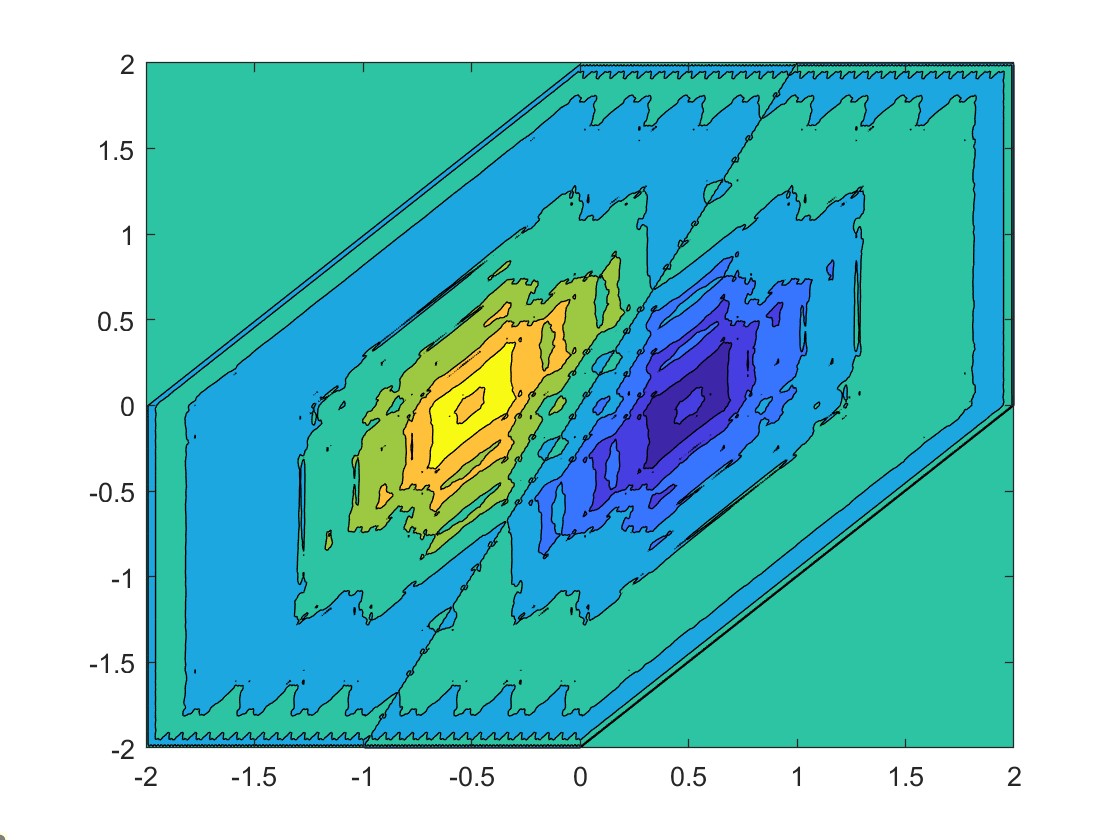}
		\caption{Contour of $\frac{\partial \phi}{\partial x}$}
	\end{subfigure}
	\caption{
		(A) is the graph of the interpolating $4I_2$-refinable function $\phi\in C^1(\R^2)$ in Example~\ref{ex1} and (B) is its contour.
		(C) is the graph of the partial derivative $\phi_{x}=\frac{\partial \phi}{\partial x} \in C^0(\R^2)$, and (D) is the contour of $\frac{\partial \phi}{\partial x}$.
	}\label{fig:ex1}
\end{figure}

\begin{exmp}\label{ex1} Let  $a_1,a_2\in l_0(\Z^2)$ be given by \er{1:a1} and \er{1:a2} where $p\in l_0(\Z^2)$ is given by \er{1:p} and $t_1,t_2\in\R$ are free parameters. Define $a\in l_0(\Z^2)$ via \er{1:a}. By imposing the $4I_2$-interpolatory constraint $a(4k)=\frac{1}{16}\td(k)$ for all $k\in\Z^2$, we obtain $t_1=\frac{t_2}{2(2t_2-1)}$. By taking $t_2=\frac{11}{64}$, the masks $a_1,a_2$ are then given by	 $$a_1=\frac{1}{672}\begin{bmatrix}0 & 0 & -11 & -22 & -11\\	0 & -22 & 106 & 106 &-22\\
	11 & 106 & 234 & 106 & -11 \\
	-22 & 106 & 106 & -22 & 0\\
	-11 & -22 & -11 & 0 & 0
	\end{bmatrix}_{[-2,2]^2},\quad a_2=	 \frac{1}{512}\begin{bmatrix}0 & 0 & 11 & 22 & 11\\
	0 & 22 & 42 & 42 & 22\\
	11 & 42 & 62 & 42 & 11\\
	22 & 42 & 42 & 22 & 0\\
	11 & 22 & 11 & 0 & 0\end{bmatrix}_{[-2,2]^2}.$$
Computation yields $\sm_2(a,4I_2)\approx 1.709055$. Using the estimates in \er{sm:a:1} with $n=1$, we obtain $\sm_\infty(a,4I_2)\ge 1.30098$. Therefore, the $2$-mask interpolatory quasi-stationary $2I_2$-subdivision scheme using the masks $a_1,a_2$ is $C^1$-convergent. See Figure~\ref{fig:ex1} for the graphs of the standard $4I_2$-refinable function $\phi$ of the mask $a$, $\frac{\partial \phi}{\partial x}$, and the contours of $\phi$ and $\frac{\partial \phi}{\partial x}$.
\end{exmp}

\subsection{Examples of $C^2$-convergent schemes}

Now we consider examples of $C^2$-convergent $r$-mask interpolatory quasi-stationary $2I_2$-subdivision schemes. We first observe that masks with $1$-ring stencil cannot achieve $C^2$-convergence. Indeed, suppose $a_1,\dots,a_r\in l_0(\Z^2)$ are $2I_2$-interpolatory and supported inside $[-2,2]^2$, then by \cite[Theorem 3.4]{hj06}, we have $\sr(a_l,2I_2)\le 2$ for all $l=1,\dots,r$. This implies that the mask $a:=2^{-2r}\sd_{a_r,2I_2}\dots\sd_{a_1,2I_2}\td$ is supported inside
$$(2^{r-1}+\dots+2+1)[-2,2]^2=[2-2^{r+1},2^{r+1}-2]^2,$$
and satisfies $\sr(a,2^rI_2)\le 2$. Hence, by \cite[Theorem 3.4]{hj06} again yields $\sm_\infty(a,2^rI_2)\le 2$ and the standard $2^rI_2$-refinable function $\phi$ of the mask $a$ is not in $C^2(\R^2)$, that is, the $r$-mask quasi-stationary $2I_2$-subdivision scheme that uses the masks $a_1,\dots,a_r$ is not $C^2$-convergent. Consequently, a $C^2$-convergent interpolatory $r$-mask quasi-stationary $2I_2$-subdivision scheme must be at least $2$-ring stencils. Here we present $C^2$-convergent examples with $r=2$.

\vspace{0.2cm}

Let $\mathcal{G}=D_4$ and parameterize two masks $a_1,a_2\in l_0(\Z^2)$ such that
\begin{itemize}
	\item $\wh{a_1}(0)=\wh{a_2}(0)=1$, $\sr(a_1,2I_2)=\sr(a_2,2I_2)=4$;
	
	\item $a_1$ and $a_2$ have $2$-ring stencil;
	
	\item $a_1$ and $a_2$ are $D_4$-symmetric about $(0,0)$;
\end{itemize}
as follows:
\begin{align}
&\wh{a_1}(\xi_1,\xi_2):=\frac{1}{256}e^{2i(\xi_1+\xi_2)}
(1+e^{-i\xi_1})^4(1+e^{-i\xi_2})^4
\wh{p_1}(\xi_1,\xi_2),
\label{D4:2:a1}\\
&\wh{a_2}(\xi_1,\xi_2):=
\frac{1}{256}e^{2i(\xi_1+\xi_2)}
(1+e^{-i\xi_1})^4(1+e^{-i\xi_2})^4
\wh{p_2}(\xi_1,\xi_2),
\label{D4:2:a2}
\end{align}	
where $t_1, \ldots, t_8\in \R$ are free parameters and $p_1,p_2\in l_0(\Z^2)$ are given by
\be\label{d4:p1}p_1=\begin{bmatrix}
t_5 &t_4 &t_3 &t_4 &t_5\\
t_4 &t_2 &t_1 &t_2 &t_4\\
t_3 &t_1 &1 - 4t_1 - 4t_5 - 4t_3 - 8t_4 - 4t_2
&t1 &t_3\\
t_4 &t_2 &t_1 &t_2 &t_4\\
t_5 &t_4 &t_3 &t_4 &t_5
\end{bmatrix}_{[-2,2]^2},\ee
and
\be\label{d4:p2}p_2=\begin{bmatrix}
0 &0 &t_8 &0 &0\\
0 &t_7 &t_6 &t_7 &t_0\\
t_8 &t_6 &1 - 4t_6 - 4t_7 - 4t_8
&t6 &t_8\\
0 &t_7 &t_6 &t_7 &0\\
0 &0 &t_8 &0 &0
\end{bmatrix}_{[-2,2]^2}.\ee
Define $a\in l_0(\Z^2)$ via \er{1:a} with $a_1,a_2\in l_0(\Z^2)$ being given by \er{D4:2:a1} and \er{D4:2:a2}. To estimate $\sm_\infty(a,4I_2)$, it suffices to estimate
$$\rho_4(a,4I_2,\mu)_\infty:=\lim_{n\to\infty}\left\|\nabla^\mu\sd_{a,4I_2}^n\td\right\|_{e_\infty(\Z^2)}^{\frac{1}{n}},\quad\text{for }\mu=(4,0),\,(3,1),\,(2,2),\,(1,3),\,(0,4).$$
Since $a_1,a_2$ are $D_4$-symmetric, we have
$$\rho_4(a,4I_2,(4,0))_\infty=\rho_4(a,4I_2,(0,4))_\infty,\quad \rho_4(a,4I_2,(3,1))_\infty=\rho_4(a,4I_2,(1,3))_\infty,$$
so
\be\label{sr4:rho}\rho_4(a,4I_2)_\infty=\max\{\rho_4(a,4I_2,(4,0))_\infty,\,\rho_4(a,4I_2,(3,1))_\infty,\,\rho_4(a,4I_2,(2,2))_\infty\}.\ee
For every $u\in l_0(\Z^2)$ and $k\in\Z^2$, note that
\be\label{diff:u}\nabla_{e_1}\nabla_{e_2}u(k)=\nabla_{e_2}u(k)-\nabla_{e_2}u(k-e_1),\quad\nabla_{e_2}^2u(k)=\nabla_{e_2}u(k)-\nabla_{e_2}u(k-e_2),\ee
from which we obtain
\be\label{est:31}\rho_4(a,4I_2,(3,1))_\infty=\lim_{n\to\infty}\left\|\nabla^2_{e_1}\nabla_{e_1}\nabla_{e_2}\sd_{a,4I_2}^n\td\right\|_{l_\infty(\Z^2)}^{\frac{1}{n}}\le\lim_{n\to\infty}\left\|\nabla^2_{e_1}\nabla_{e_2}\sd_{a,4I_2}^n\td\right\|_{l_\infty(\Z^2)}^{\frac{1}{n}},\ee
\be\label{est:22}\rho_4(a,4I_2,(2,2))_\infty=\lim_{n\to\infty}\left\|\nabla^2_{e_1}\nabla_{e_2}^2\sd_{a,4I_2}^n\td\right\|_{l_\infty(\Z^2)}^{\frac{1}{n}}\le \lim_{n\to\infty}\left\|\nabla^2_{e_1}\nabla_{e_2}\sd_{a,4I_2}^n\td\right\|_{l_\infty(\Z^2)}^{\frac{1}{n}}.\ee
Hence
$$\rho_4(a,4I_2)_\infty\le\max\left\{\rho_4(a,4I_2,(4,0))_\infty,\,\lim_{n\to\infty}\left\|\nabla^2_{e_1}\nabla_{e_2}\sd_{a,4I_2}^n\td\right\|_{l_\infty(\Z^2)}^{\frac{1}{n}}\right\}.$$
Let $p_1,p_2\in l_0(\Z^2)$ be given by \er{d4:p1} and \er{d4:p2}. Define $h_1,h_2\in l_0(\Z^2)$ via
\be\label{h1}\wh{h_1}(\xi_1,\xi_2):=\frac{1}{65536}e^{6i(\xi_1+\xi_2)}(1+e^{-i\xi_2})^4(1+e^{-2i\xi_2})^4\wh{p_1}(2\xi_1,2\xi_2)\wh{p_2}(\xi_1,\xi_2),\ee
\be\label{h2}\wh{h_2}(\xi_1,\xi_2):=\frac{1}{65536}e^{6i(\xi_1+\xi_2)}(1+e^{-i\xi_1})^2(1+e^{-2i\xi_1})^2(1+e^{-i\xi_2})^3(1+e^{-2i\xi_2})^3\wh{p_1}(2\xi_1,2\xi_2)\wh{p_2}(\xi_1,\xi_2).\ee
Note that $$\wh{a}(\xi_1,\xi_2)=\frac{\wh{\nabla^{(4,0)}\td}(4\xi_1,4\xi_2)}{\wh{\nabla^{(4,0)}\td}(\xi_1,\xi_2)}\wh{h_1}(\xi_1,\xi_2)=\frac{\wh{\nabla^2_{e_1}\nabla_{e_2}\td}(4\xi_1,4\xi_2)}{\wh{\nabla^2_{e_1}\nabla_{e_2}\td}(\xi_1,\xi_2)}\wh{h_2}(\xi_1,\xi_2).$$
We then conclude from Theorems~\ref{thm:sm:1} and ~\ref{thm:sm:2} that
$$\rho_4(a,4I_2,(4,0))_\infty=\lim_{n\to\infty}\left\|\sd_{h_1,4I_2}^n\td\right\|_{l_\infty(\Z^2)}^{\frac{1}{n}}\le \left(\sup_{\gamma\in\Gamma_{4^n}}\sum_{k\in\Z^2}\left|\sd_{h_1,4I_2}^n\td(\gamma+4^nk)\right|\right)^{\frac{1}{n}},$$
$$\lim_{n\to\infty}\left\|\nabla^2_{e_1}\nabla_{e_2}\sd_{a,4I_2}^n\td\right\|_{l_\infty(\Z^2)}^{\frac{1}{n}}=\lim_{n\to\infty}\left\|\sd_{h_2,4I_2}^n\td\right\|_{l_\infty(\Z^2)}^{\frac{1}{n}}\le \left(\sup_{\gamma\in\Gamma_{4^n}}\sum_{k\in\Z^2}\left|\sd_{h_2,4I_2}^n\td(\gamma+4^nk)\right|\right)^{\frac{1}{n}},$$
for all $n\in\N$. Consequently, we have
\be\label{est:sm:2}\sm_\infty(a,4I_2)\ge-\log_4\left(\max_{j=1,2}\left(\sup_{\gamma\in\Gamma_{4^n}}\sum_{k\in\Z^2}\left|\sd_{h_j,4I_2}^n \td(\gamma+4^nk)\right|\right)^{\frac{1}{n}}\right),\quad\forall n\in\N,
\ee
where $h_1,h_2\in l_0(\Z^2)$ are given by \er{h1} and \er{h2}.

\begin{exmp}\label{ex4}Let $p_1,p_2\in l_0(\Z^2)$ be given by \er{d4:p1} and \er{d4:p2} for some free parameters $t_1,\dots,t_8\in\R$. Let $a_1,a_2\in l_0(\Z^2)$ be given by \er{D4:2:a1} and \er{D4:2:a2}. Define $a\in l_0(\Z^2)$ via \er{1:a}. By imposing the $4I_2$-interpolatory constraint $a(4k)=\frac{1}{16}\td(k)$ for all $k\in\Z^2$, we have many solutions and here we present two solutions. The first choice is $t_3=t_4=t_5 = 0$, $t_8 = t$, and
\begin{align*}
&t_1=-\tfrac{90035}{1027628}t^4 - \tfrac{21065363}{12331536}t^3 + \tfrac{10201651}{3082884}t^2 +\tfrac{13085857}{6165768}t - \tfrac{640565}{1541442},\\
&t_2 = -\tfrac{279697}{4110512}t^4 - \tfrac{65258017}{49326144}t^3 + \tfrac{16876337}{6165768}t^2 + \tfrac{44294515}{12331536}t + \tfrac{2660345}{3082884},\\
%&t_3=t_4=t_5 = 0,\\
&t_6 = \tfrac{1398881}{513814}t^4 + \tfrac{325752977}{6165768}t^3 - \tfrac{84880928}{770721}t^2
- \tfrac{85869100}{770721}t - \tfrac{15771775}{770721},\\
&t_7 = -\tfrac{469557}{513814}t^4 - \tfrac{36394039}{2055256}t^3 + \tfrac{28901470}{770721}t^2
+ \tfrac{9016005}{256907}t + \tfrac{19610825}{3082884},%\quad t_8 = t,
\end{align*}
where $t\approx -0.2395777$ is a root of $132t^5 + 2651t^4 - 3600t^3 - 8896t^2 - 4560t - 640=0$. We have
$$t_1 \approx-0.71089031,\quad t_2 \approx 0.17745551,\quad t_6 \approx -0.810183, \quad t_7\approx 0.3462063,\quad  t_8 \approx-0.23957771.$$
The two masks $a_1,a_2$ are then approximately given by
\begin{align*}
&a_1={\fontsize{9}{9}\selectfont \begingroup % keep the change local
\setlength\arraycolsep{1pt}\begin{bmatrix}0.0006929244 & -4.6872\times 10^{-6} & -0.0062550684 & -0.0111149136 & -0.0062550684 & -4.6872 \times 10^{-6} & 0.0006929244\\
-4.6872\times 10^{-6} & 0.0011210220 & 0.0045262728 & 0.0068011272 & 0.0045262728 & 0.0011210220 & -4.6872\times 10^{-6}\\
-0.0062550684 & 0.0045262728 & 0.0744007068 & 0.1272387312 & 0.0744007068 & 0.0045262728 & -0.0062550684\\
-0.0111149136 & 0.0068011272 & 0.1272387312 & 0.2186453808 & 0.1272387312 & 0.0068011272 & -0.0111149136\\
-0.0062550684 & 0.0045262728 & 0.0744007068 & 0.1272387312 & 0.0744007068 & 0.0045262728 & -0.0062550684\\
-4.6872\times 10^{-6} & 0.0011210220 & 0.0045262728 & 0.0068011272 & 0.0045262728 & 0.0011210220 & -4.6872\times 10^{-6}\\
0.0006929244 & -4.6872\times 10^{-6} & -0.0062550684 & -0.0111149136 & -0.0062550684 & -4.6872\times 10^{-6} & 0.0006929244\end{bmatrix}_{[-3,3]^2},\endgroup}\end{align*}
\begin{align*}
&a_2={\fontsize{9}{9}\selectfont \begingroup % keep the change local
\setlength\arraycolsep{1pt}\begin{bmatrix}0 & 0 & -0.00093840 & -0.00375360 & -0.00563040 & -0.00375360 & -0.00093840 & 0 & 0\\
0 & 0.0013294 & -0.00234600 & -0.02134860 & -0.03534640 & -0.02134860 & -0.00234600 & 0.0013294 & 0\\
-0.00093840 & -0.00234600 & -0.00602140 & -0.01798600 & -0.02674440 & -0.01798600 & -0.00602140 & -0.00234600 & -0.00093840\\
-0.00375360 & -0.02134860 & -0.01798600 & 0.06013580 & 0.12105360 & 0.06013580 & -0.01798600 & -0.02134860 & -0.00375360\\
-0.00563040 & -0.03534640 & -0.02674440 & 0.12105360 & 0.23616400 & 0.12105360 & -0.02674440 & -0.03534640 & -0.00563040\\
-0.00375360 & -0.02134860 & -0.01798600 & 0.06013580 & 0.12105360 & 0.06013580 & -0.01798600 & -0.02134860 & -0.00375360\\
-0.00093840 & -0.00234600 & -0.00602140 & -0.01798600 & -0.02674440 & -0.01798600 & -0.00602140 & -0.00234600 & -0.00093840\\
0 & 0.0013294 & -0.00234600 & -0.02134860 & -0.03534640 & -0.02134860 & -0.00234600 & 0.0013294 & 0\\
0 & 0 & -0.00093840 & -0.00375360 & -0.00563040 & -0.00375360 & -0.00093840 & 0 & 0\end{bmatrix}_{[-4,4]^2}.\endgroup}\end{align*}
Moreover, direct computation yields $\sm_2(a,4I_2)\approx 2.616519$. Using the estimation \er{est:sm:2} with $n=2$, we obtain $\sm_\infty(a,4I_2)\ge 2.07607$. Therefore, the $2$-mask interpolatory quasi-stationary $2I_2$-subdivision scheme using the above masks $a_1,a_2$ is $C^2$-convergent.

%\vspace{0.2cm}

Another choice is $t_4=0, t_5=-\tfrac{5}{128}, t_7 = \tfrac{t}{4}, t_8=0$ and
\begin{align*}
%&t_4=0,\quad t_5=-\tfrac{5}{128},\quad t_7 = \tfrac{t}{4},\quad t_8=0,\\
&t_1 = \tfrac{424994944920607}{24908728622815420334080}t^6 + \tfrac{7975356264181027}{9963491449126168133632}t^5 - \tfrac{143174662095446521}{9963491449126168133632}t^4 -\tfrac{201920022567847971491}{24908728622815420334080}t^3 \\
&\qquad - \tfrac{1877989942288756509}{6227182155703855083520}
t^2 -
\tfrac{6489286410377804222101}{9963491449126168133632}t -\tfrac{26798558462517617887845}{9963491449126168133632},\\
&t_2 = \tfrac{113900026967}{77839776946298188544}t^6 + \tfrac{6721822856723}{155679553892596377088}t^5 - \tfrac{362570802099563}{155679553892596377088}t^4
- \tfrac{52158296539168801}{77839776946298188544}t^3\\
&\qquad  + \tfrac{464437583880593005}{38919888473149094272}t^2 - \tfrac{15627987102075198505}{155679553892596377088}t
+\tfrac{203270998898608942625}{155679553892596377088},\\
&t_3 = -\tfrac{2633158952204571}{498174572456308406681600}t^6 -\tfrac{46155986031448111}{199269828982523362672640}t^5 +\tfrac{1033152943943700317}{199269828982523362672640}t^4 +\tfrac{1243781358911040235183}{498174572456308406681600}
t^3\\
&\qquad - \tfrac{949630637866665816943}{124543643114077101670400}t^2
+\tfrac{43309282197078887566441}{199269828982523362672640}t -\tfrac{3476955096904628746267}{39853965796504672534528},\\
&t_6 = \tfrac{561344301981}{121624651478590919600}t^6 + \tfrac{7723586258751}{48649860591436367840}t^5 - \tfrac{166473281489741}{24324930295718183920}t^4
- \tfrac{262791785728360213}{121624651478590919600}t^3 \\
&\qquad + \tfrac{3373757083893137367}{121624651478590919600}t^2
- \tfrac{27687970920178103501}{48649860591436367840}t + \tfrac{2232880122164793821}{4864986059143636784},
\end{align*}
where $t\approx 2.233641927$ is a root of
\[
2t^7 + 95t^6 - 1660t^5 - 952671t^4 - 573006t^3 - 61196575t^2 - 340415800t + 1095865375=0.
\]
We have
$$t_1\approx-4.2366142, \quad t_2\approx1.1334896,\quad t_3\approx 0.38811383, \quad t_6\approx-0.69810232,\quad t_7\approx0.55841048,$$
and the masks $a_1, a_2$ are approximately given by
\begin{align*}
&a_1={\fontsize{9}{9}\selectfont \begingroup % keep the change local
\setlength\arraycolsep{1pt}\begin{bmatrix}-0.0001521 & -0.0006084 & 0.0005694 & 0.0053196 & 0.0085878 & 0.0053196 & 0.0005694 & -0.0006084 & -0.0001521\\
-0.0006084 & 0.0018564 & 0.0030576 & -0.0142116 & -0.0296088 & -0.0142116 & 0.0030576 & 0.0018564 & -0.0006084\\
0.0005694 & 0.0030576 & -0.0032916 & -0.0315744 & -0.0515892 & -0.0315744 & -0.0032916 & 0.0030576 & 0.0005694\\
0.0053196 & -0.0142116 & -0.0315744 & 0.0867204 & 0.1975272 & 0.0867204 & -0.0315744 & -0.0142116 & 0.0053196\\
0.0085878 & -0.0296088 & -0.0515892 & 0.1975272 & 0.4218396 & 0.1975272 & -0.0515892 & -0.0296088 & 0.0085878\\
0.0053196 & -0.0142116 & -0.0315744 & 0.0867204 & 0.1975272 & 0.0867204 & -0.0315744 & -0.0142116 & 0.0053196\\
0.0005694 & 0.0030576 & -0.0032916 & -0.0315744 & -0.0515892 & -0.0315744 & -0.0032916 & 0.0030576 & 0.0005694\\
-0.0006084 & 0.0018564 & 0.0030576 & -0.0142116 & -0.0296088 & -0.0142116 & 0.0030576 & 0.0018564 & -0.0006084\\
-0.0001521 & -0.0006084 & 0.0005694 & 0.0053196 & 0.0085878 & 0.0053196 & 0.0005694 & -0.0006084 & -0.0001521\end{bmatrix}_{[-4,4]^2},\endgroup}\end{align*}
\begin{align*}
&a_2={\fontsize{9}{9}\selectfont \begingroup % keep the change local
\setlength\arraycolsep{1pt}\begin{bmatrix}0.002145 & 0.005694 & 0.003471 & -0.000156 & 0.003471 & 0.005694 & 0.002145\\
0.005694 & 0.017862 & 0.020202 & 0.016068 & 0.020202 & 0.017862 & 0.005694\\
0.003471 & 0.020202 & 0.049569 & 0.065676 & 0.049569 & 0.020202 & 0.003471\\
-0.000156 & 0.016068 & 0.065676 & 0.098904 & 0.065676 & 0.016068 & -0.000156\\
0.003471 & 0.020202 & 0.049569 & 0.065676 & 0.049569 & 0.020202 & 0.003471\\
0.005694 & 0.017862 & 0.020202 & 0.016068 & 0.020202 & 0.017862 & 0.005694\\
0.002145 & 0.005694 & 0.003471 & -0.000156 & 0.003471 & 0.005694 & 0.002145\end{bmatrix}_{[-3,3]^2}.\endgroup}\end{align*}
Moreover, direct computation yields $\sm_2(a,4I_2)\approx 3.074404$. Hence, using \er{sm:inf:2} we obtain $\sm_\infty(a,4I_2)\ge 2.074404$. Therefore, the $2$-mask interpolatory quasi-stationary $2I_2$-subdivision scheme using the above masks $a_1,a_2$ is $C^2$-convergent. See Figure~\ref{fig:ex4} for the graphs of the $4I_2$-standard refinable function $\phi$ of the mask $a$ with the second choice, $\frac{\partial^2 \phi}{\partial x^2}$, and the contours of $\phi$ and $\frac{\partial^2 \phi}{\partial x^2}$.

\begin{figure}[htbp]
\centering
\begin{subfigure}[b]{0.24\textwidth} \includegraphics[width=\textwidth,height=0.8\textwidth]{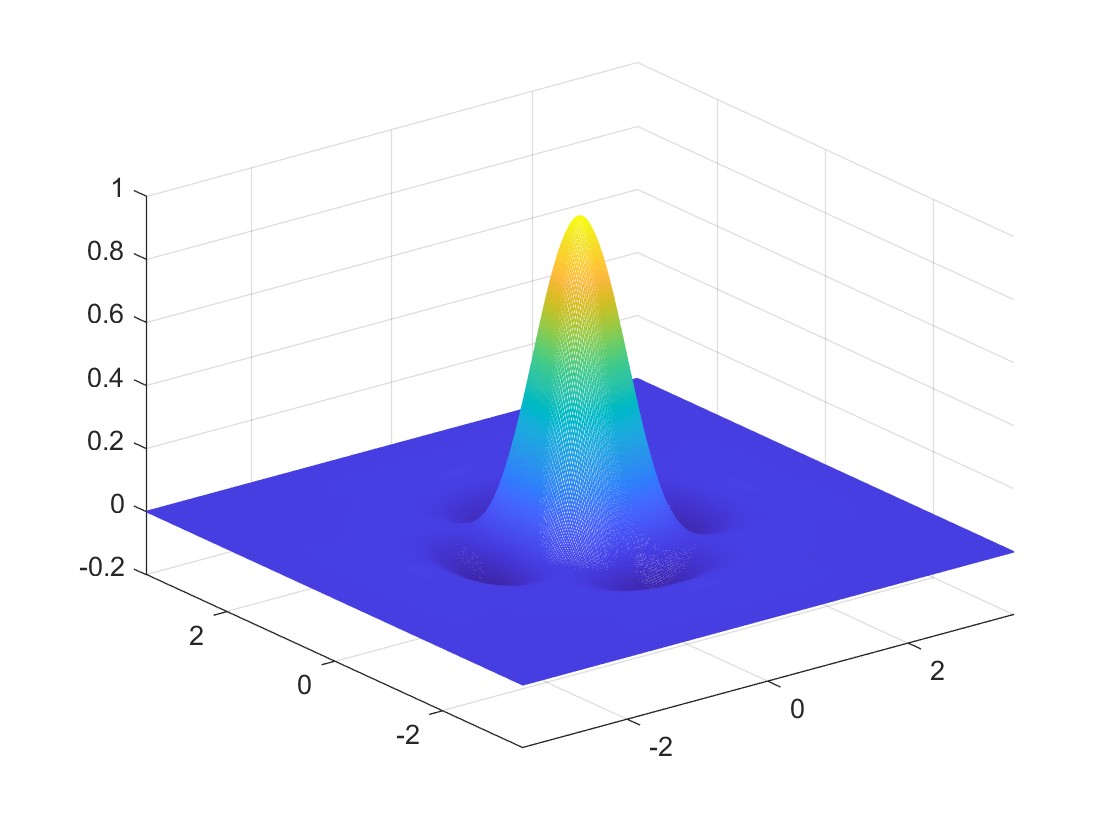}
	\caption{$\phi\in \mathscr{C}^2(\R^2)$}
\end{subfigure}
\begin{subfigure}[b]{0.22\textwidth} \includegraphics[width=\textwidth,height=0.8\textwidth]{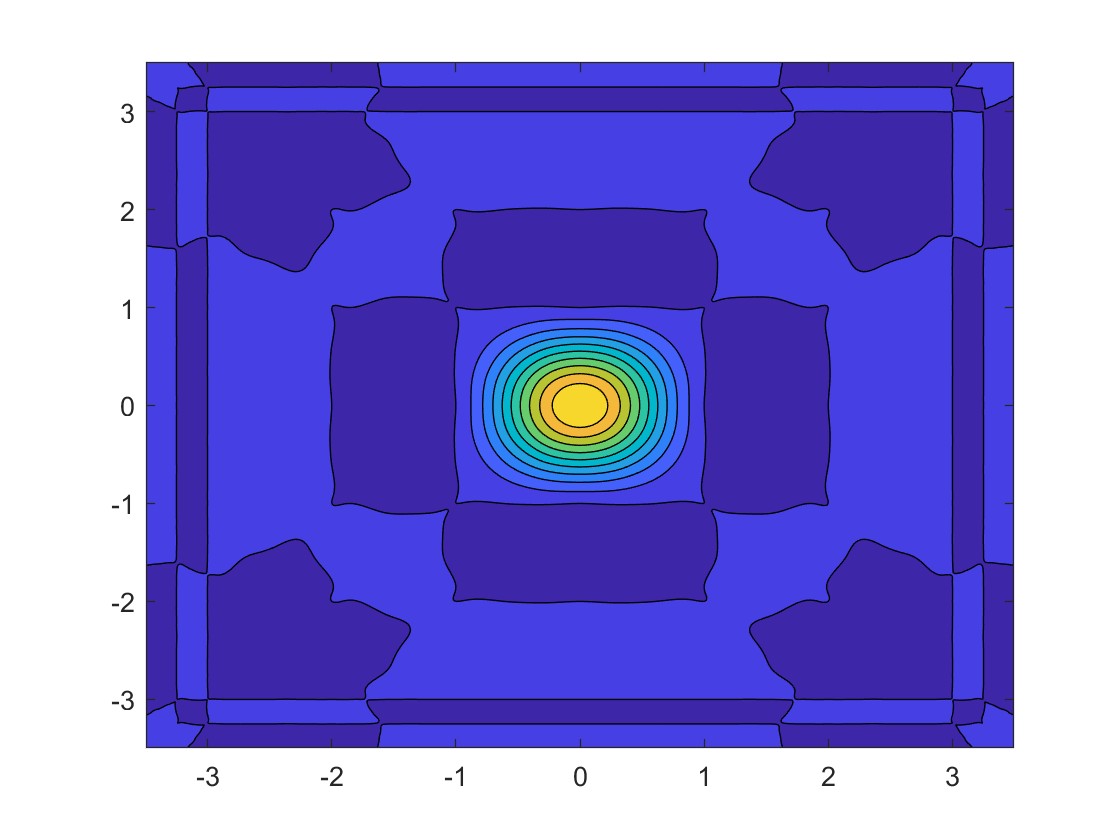}
	\caption{Contour of $\phi$}
\end{subfigure}
\begin{subfigure}[b]{0.24\textwidth} \includegraphics[width=\textwidth,height=0.8\textwidth]{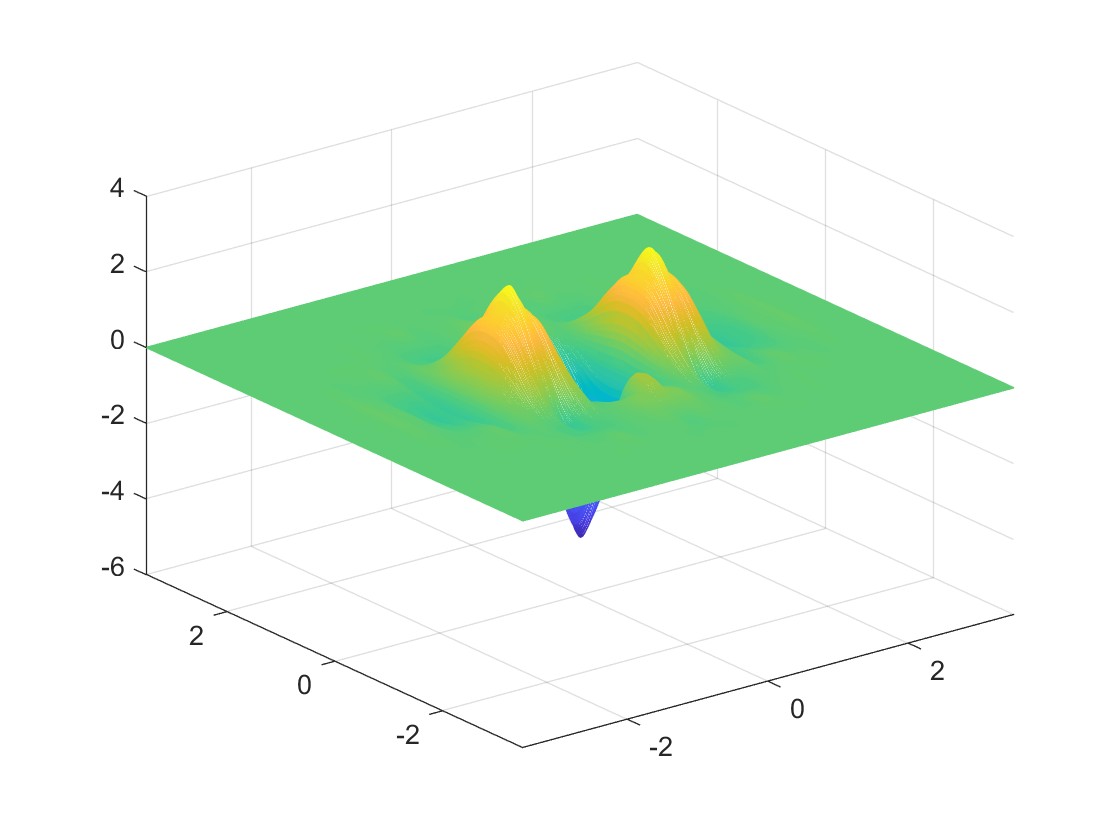} \caption{$\frac{\partial^2 \phi}{\partial x^2} \in C^0(\R^2)$}
\end{subfigure}
\begin{subfigure}[b]{0.22\textwidth} \includegraphics[width=\textwidth,height=0.8\textwidth]{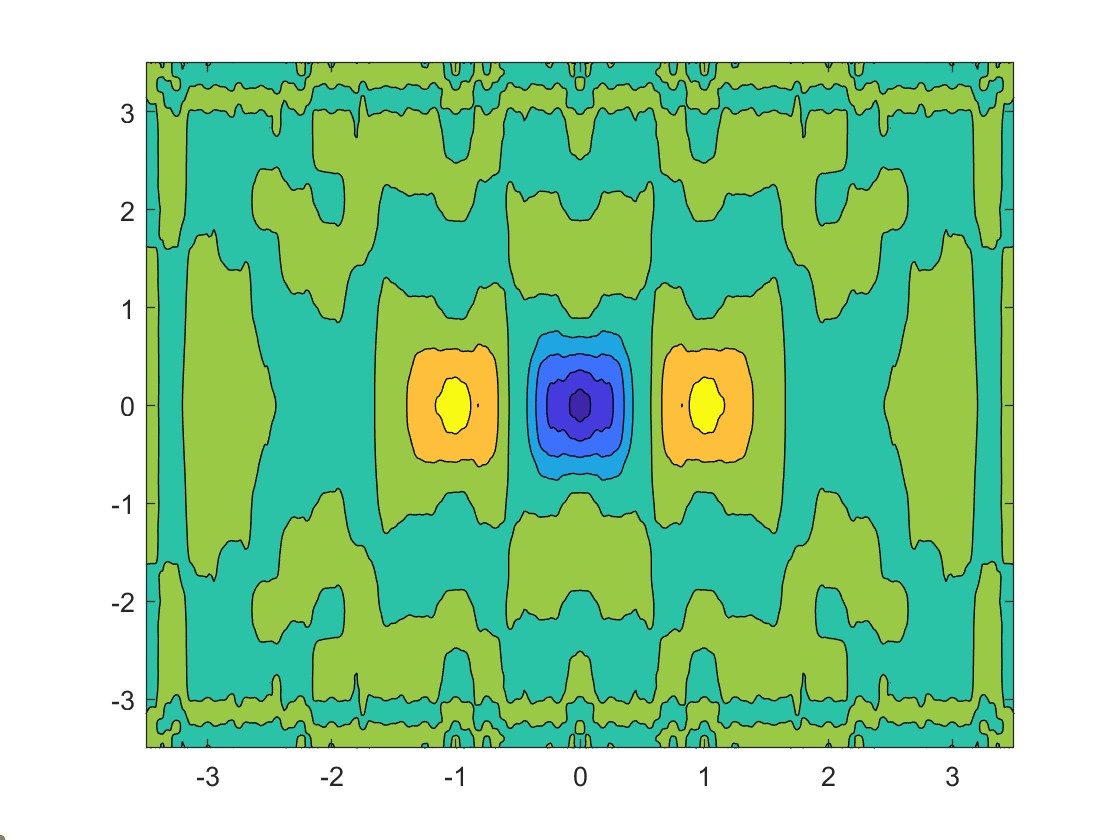}
	\caption{Contour of $\frac{\partial^2 \phi}{\partial x^2}$}
\end{subfigure}
\caption{
	(A) is the graph of the interpolating $4I_2$-refinable function $\phi\in C^2(\R^2)$ in Example~\ref{ex4} and (B) is its contour.
	(C) is the graph of the partial derivative $\frac{\partial^2 \phi}{\partial x^2} \in C^0(\R^2)$, and (D) is the contour of $\frac{\partial^2 \phi}{\partial x^2}$.
}\label{fig:ex4}
\end{figure}

\end{exmp}

Let $\mathcal{G}=D_6$ and parameterize two masks $a_1,a_2\in l_0(\Z^2)$ such that  \begin{itemize}
	\item $\wh{a_1}(0)=\wh{a_2}(0)=1$, $\sr(a_1,2I_2)=\sr(a_2,2I_2)=4$;
	
	\item $a_1$ and $a_2$ have two-ring stencils;
	
	\item $a_1$ and $a_2$ are $D_6$-symmetric about $(0,0)$;
\end{itemize}
as follows:
\be\label{2:a1}\wh{a_1}(\xi_1,\xi_2):=\frac{1}{64}(1+e^{-i\xi_1})^2(1+e^{-i\xi_2})^2(1+e^{i(\xi_1+\xi_2)})^2\wh{g_1}(\xi_1,\xi_2),\ee
\be\label{2:a2}\wh{a_2}(\xi_1,\xi_2):=\frac{1}{64}(1+e^{-i\xi_1})^2(1+e^{-i\xi_2})^2(1+e^{i(\xi_1+\xi_2)})^2\wh{g_2}(\xi_1,\xi_2),\ee
where $g_1,g_2\in l_0(\Z^2)$ are given by
\be\label{g1:g2}g_1=\begin{bmatrix}0 & 0 & t_3 & t_2 & t_3\\
	0 & t_2 & t_1 & t_1 & t_2\\
	t_3 & t_1 & 1-6t_3-6t_2-6t_1 & t_1 & t_3\\
	t_2 & t_1 & t_1 & t_2 & 0\\
	t_3 & t_2 & t_3 & 0 & 0\end{bmatrix}_{[-2,2]^2},\quad g_2=\begin{bmatrix}0 & 0 & t_6 & t_5 & t_6\\
	0 & t_5 & t_4 & t_4 & t_5\\
	t_6 & t_4 & 1-6t_6-6t_5-6t_4 & t_4 & t_6\\
	t_5 & t_4 & t_4 & t_5 & 0\\
	t_6 & t_5 & t_6 & 0 & 0\end{bmatrix}_{[-2,2]^2},\ee
where $t_1,\dots,t_6\in\R$ are free parameters. Define $a\in l_0(\Z^2)$ via \er{1:a}. Since $a_1,a_2$ are $D_6$-symmetric about $(0,0)$, it is easy to see that \er{sr4:rho} holds. By letting $E=\begin{bmatrix}1 &-1\\ 0 &-1\end{bmatrix}\in D_6$, \er{e1:e2} holds, which together with \er{n:e1}, yields $\rho_4(a,4I_2,(4,0))_\infty\le \rho_4(a,4I_2,(3,1))_\infty.$ Furthermore, one can conclude from \er{diff:u} that \er{est:31} and \er{est:22} must hold. Therefore, we have
$$\rho_4(a,4I_2)_\infty\le \lim_{n\to\infty}\left\|\nabla^2_{e_1}\nabla_{e_2}\sd_{a,4I_2}^n\td\right\|_{l_\infty(\Z^2)}^{\frac{1}{n}}.$$
Let $g_1,g_2\in\l_0(\Z^2)$ be given by \er{g1:g2}. Define $h\in l_0(\Z^2)$ via
\be\label{hh}\wh{h}(\xi_1,\xi_2)=\frac{1}{4096}(1+e^{-i\xi_2})(1+e^{-2i\xi_2})(1+e^{i(\xi_1+\xi_2)})^2(1+e^{2i(\xi_1+\xi_2)})^2\wh{g_1}(2\xi_1,2\xi_2)\wh{g_2}(\xi_1,\xi_2).\ee
Note that
$$\wh{a}(\xi_1,\xi_2)=\frac{\wh{\nabla^2_{e_1}\nabla_{e_2}\td}(4\xi_1,4\xi_2)}{\wh{\nabla^2_{e_1}\nabla_{e_2}\td}(\xi_1,\xi_2)}\wh{h}(\xi_1,\xi_2).$$
It then follows from Theorems~\ref{thm:sm:1} and ~\ref{thm:sm:2} that
$$\lim_{n\to\infty}\left\|\nabla^2_{e_1}\nabla_{e_2}\sd_{a,4I_2}^n\td\right\|_{l_\infty(\Z^2)}^{\frac{1}{n}}=\lim_{n\to\infty}\left\|\sd_{h,4I_2}^n\td\right\|_{l_\infty(\Z^2)}^{\frac{1}{n}}\le\left(\sup_{\gamma\in\Gamma_{4^n}}\sum_{k\in\Z^2}\left|\sd_{h,4I_2}^n\td(\gamma+4^nk)\right|\right)^{\frac{1}{n}},\quad\forall n\in\N.$$
Consequently, we have
\be\label{sm:a:2}\sm_\infty(a,4I_2)\ge -\log_4\left(\sup_{\gamma\in\Gamma_{4^n}}\sum_{k\in\Z^2}\left|\sd_{h,4I_2}^n\td(\gamma+4^nk)\right|\right)^{\frac{1}{n}},\quad\forall n\in\N,\ee
where $h\in l_0(\Z^2)$ is given by \er{hh}.

\begin{exmp}\label{ex2} Let $g_1,g_2\in l_0(\Z^2)$ be given by \er{g1:g2} where $t_1,\dots,t_6\in\R$ are free parameters. Define $a_1,a_2\in l_0(\Z^2)$ via \er{2:a1} and \er{2:a2} and define $a\in l_0(\Z^2)$ via \er{1:a}. By imposing the $4I_2$-interpolatory constraint $a(4k)=\frac{1}{16}\td(k)$ for all $k\in\Z^2$, we have many solutions and here we present two solutions. The first choice is $$t_1=-\frac{5}{8},\quad t_2=\frac{5}{48},\quad t_3=t_4=t_5=t_6=0.$$
The two masks $a_1,a_2$ are then given by
	$$a_1=\frac{1}{3072}\begin{bmatrix}0 & 0 & 0 & 0 & 0 & 5 & 10 & 5 & 0\\
0 & 0 & 0 & 5 & -10 & -55 & -55 & -10 & 5\\
0 & 0 & 10 & -55 & -42 & 46 & -42 & -55 & 10\\
0 & 5 & -55 & 46 & 448 & 448 & 46 & -55 & 5\\
0 & -10 & -42 & 448 & 960 & 448 & -42 & -10 & 0\\
5 & -55 & 46 & 448 & 448 & 46 & -55 & 5 & 0\\
10 & -55 & -42 & 46 & -42 & -55 & 10 & 0 & 0\\
5 & -10 & -55 & -55 & -10 & 5 & 0 & 0 & 0\\
0 & 5 & 10 & 5 & 0 & 0 & 0 & 0 & 0\end{bmatrix}_{[-4,4]^2},$$
$$a_2=\frac{1}{64}\begin{bmatrix}0 & 0 & 1 & 2 & 1\\
0 & 2 & 6 & 6 & 2\\
1 & 6 & 10 & 6 & 1\\
2 & 6 & 6 & 2 & 0\\
1 & 2 & 1 & 0 & 0\end{bmatrix}_{[-2,2]^2}.$$
Computaiton yields $\sm_2(a,4I_2)\approx2.653820$. Using the estimation \er{sm:a:2} with $n=2$, we obtain $\sm_\infty(a,4I_2)\ge 2.06210$. Therefore, the $2$-mask interpolatory quasi-stationary $2$-subdivision scheme using the above masks $a_1,a_2$ is $C^2$-convergent. See the first row of Figure~\ref{fig:ex2} for the graphs of the standard $4I_2$-refinable function $\phi$ of the mask $a$, $\frac{\partial^2 \phi}{\partial x^2}$ and the contours of $\phi$ and $\phi_{xx}$.\\

%\vspace{0.2cm}
	
Another choice is
\[
t_1=-\tfrac{2t^2}{7}-\tfrac{145t}{112}-\tfrac{151}{112},\quad t_2=\tfrac{2t^2}{21}+\tfrac{131t}{336}+\tfrac{109}{336},\quad t_3=0,\quad t_4=\tfrac{t}{2},\quad t_5=\tfrac{1}{4},\quad t_6=-\tfrac{1}{8},
\]
where $t\approx -0.133008$ is a root of $32t^3 + 141t^2 + 146t + 17=0$. The two masks $a_1,a_2$ are then approximately given by
\begin{align*}
&a_1={\fontsize{10}{10}\selectfont \begingroup % keep the change local
	\setlength\arraycolsep{1pt}\begin{bmatrix}0 & 0 & 0 & 0 & 0 & 0.0042744 & 0.0085488 & 0.0042744 & 0\\
	0 & 0 & 0 & 0.0042744 & -0.0013104 & -0.0253032 & -0.0253032 & -0.0013104 & 0.0042744\\
	0 & 0 & 0.0085488 & -0.0253032 & -0.0323232 & 0.0030576 & -0.0323232 & -0.0253032 & 0.0085488\\
	0 & 0.0042744 & -0.0253032 & 0.0030576 & 0.1656096 & 0.1656096 & 0.0030576 & -0.0253032 & 0.0042744\\
	0 & -0.0013104 & -0.0323232 & 0.1656096 & 0.3932448 & 0.1656096 & -0.0323232 & -0.0013104 & 0\\	 0.0042744 & -0.0253032 & 0.0030576 & 0.1656096 & 0.1656096 & 0.0030576 & -0.0253032 & 0.0042744 & 0\\
	0.0085488 & -0.0253032 & -0.0323232 & 0.0030576 & -0.0323232 & -0.0253032 & 0.0085488 & 0 & 0\\
	0.0042744 & -0.0013104 & -0.0253032 & -0.0253032 & -0.0013104 & 0.0042744 & 0 & 0 & 0\\
	0 & 0.0042744 & 0.0085488 & 0.0042744 & 0 & 0 & 0 & 0 & 0\end{bmatrix}_{[-4,4]^2},\endgroup}\end{align*}

\begin{align*}
&a_2={\fontsize{10}{10}\selectfont \begingroup % keep the change local
	 \setlength\arraycolsep{1pt}\begin{bmatrix}0,0 & 0 & 0 & -0.0019500 & 0 & 0.0039000 & 0 & -0.0019500\\
	0 & 0 & 0 & 0 & 0.00286260 & 0.00858780 & 0.00858780 & 0.00286260 & 0\\
	0 & 0 & 0.0039000 & 0.00858780 & 0.02315040 & 0.03692520 & 0.02315040 & 0.00858780 & 0.0039000\\
	0 & 0 & 0.00858780 & 0.03692520 & 0.06783660 & 0.06783660 & 0.03692520 & 0.00858780 & 0\\
	-0.0019500 & 0.00286260 & 0.02315040 & 0.06783660 & 0.09899760 & 0.06783660 & 0.02315040 & 0.00286260 & -0.0019500\\
	0 & 0.00858780 & 0.03692520 & 0.06783660 & 0.06783660 & 0.03692520 & 0.00858780 & 0 & 0\\
	0.0039000 & 0.00858780 & 0.02315040 & 0.03692520 & 0.02315040 & 0.00858780 & 0.0039000 & 0 & 0\\
	0 & 0.00286260 & 0.00858780 & 0.00858780 & 0.00286260 & 0 & 0 & 0 & 0\\
	-0.0019500 & 0 & 0.0039000 & 0 & -0.0019500 & 0 & 0 & 0 & 0\end{bmatrix}_{[-4,4]^2}.\endgroup}\end{align*}
Moreover, computation yields $\sm_2(a,4I_2)\approx 3.041495$ and the relation \er{sm:inf:2} yields $\sm_\infty(a,4I_2)\ge 2.041495$. Therefore, the $2$-mask interpolatory quasi-stationary $2I_2$-subdivision scheme using the above masks $a_1,a_2$ is $C^2$-convergent.
See the second row of Figure~\ref{fig:ex2} for the graphs of the standard $4I_2$-refinable function $\phi$ of the mask $a$, $\frac{\partial^2 \phi}{\partial x^2}$,  and the contours of $\phi$ and $\frac{\partial^2 \phi}{\partial x^2}$.

\begin{figure}[htbp]
	\centering
	\begin{subfigure}[b]{0.24\textwidth} \includegraphics[width=\textwidth,height=0.8\textwidth]{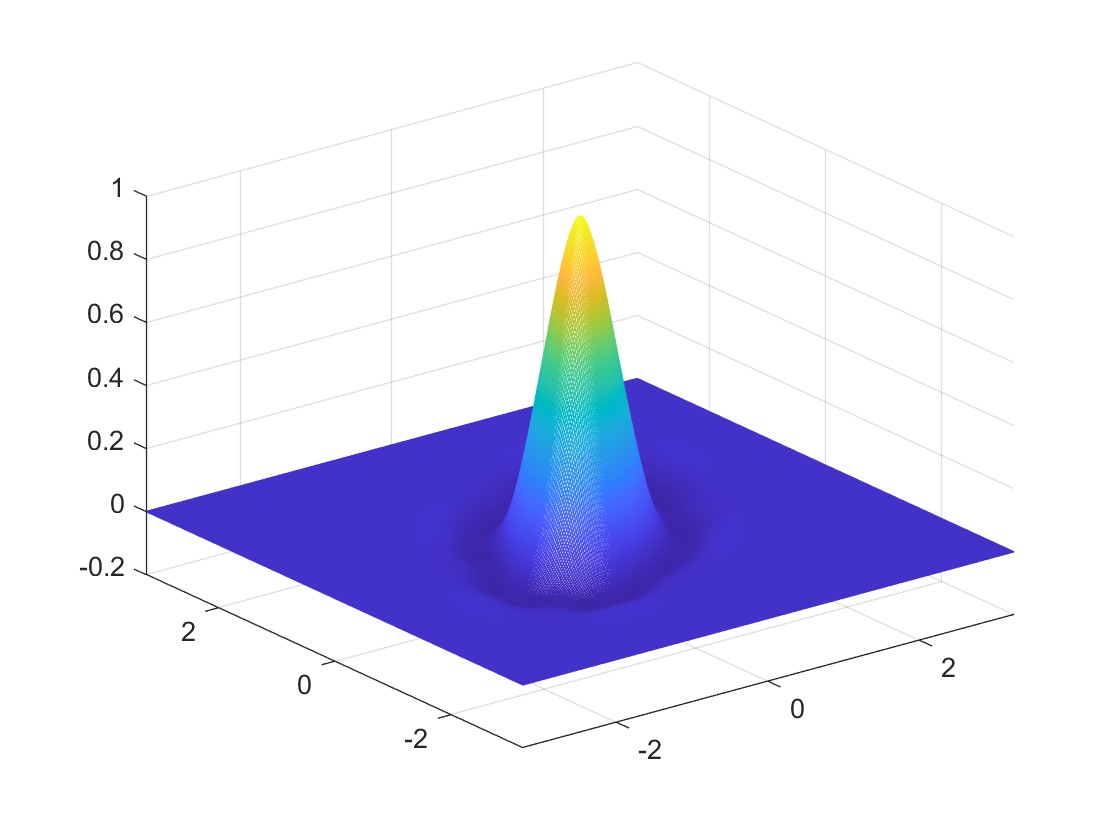}
		\caption{$\phi\in \mathscr{C}^2(\R^2)$}
	\end{subfigure}
	\begin{subfigure}[b]{0.22\textwidth} \includegraphics[width=\textwidth,height=0.8\textwidth]{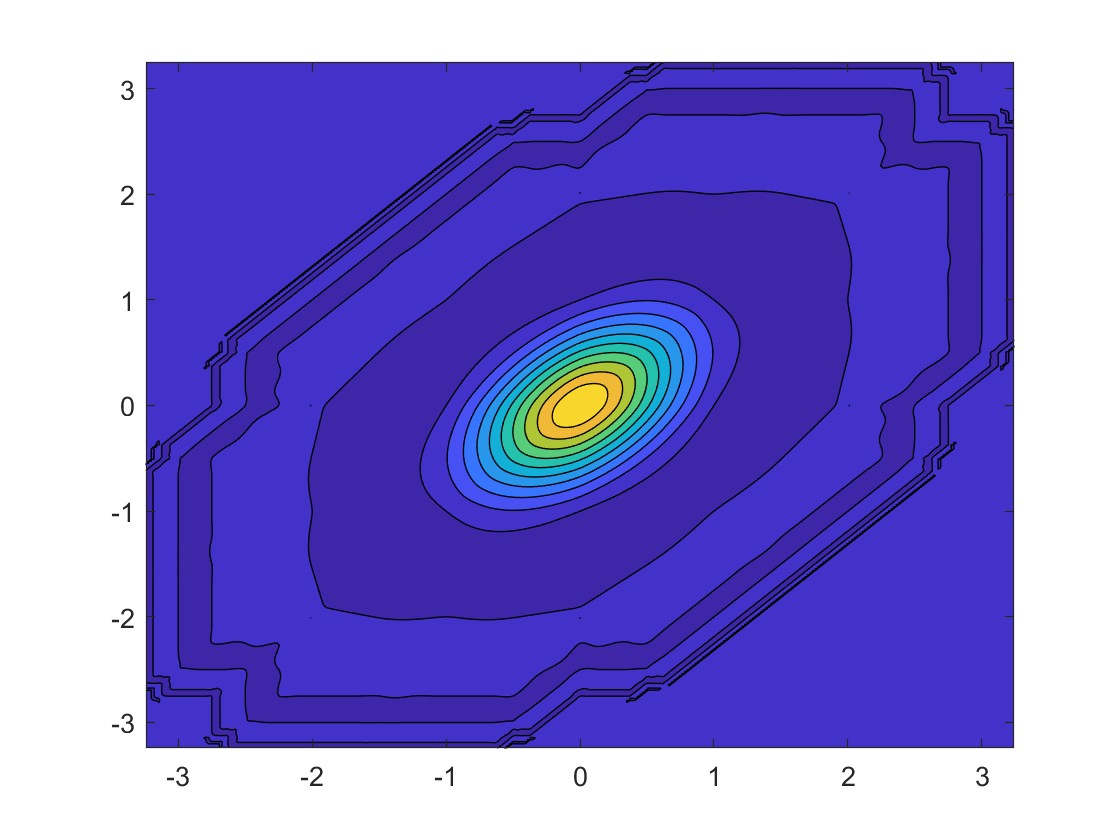}
		\caption{Contour of $\phi$}
	\end{subfigure}
	\begin{subfigure}[b]{0.24\textwidth} \includegraphics[width=\textwidth,height=0.8\textwidth]{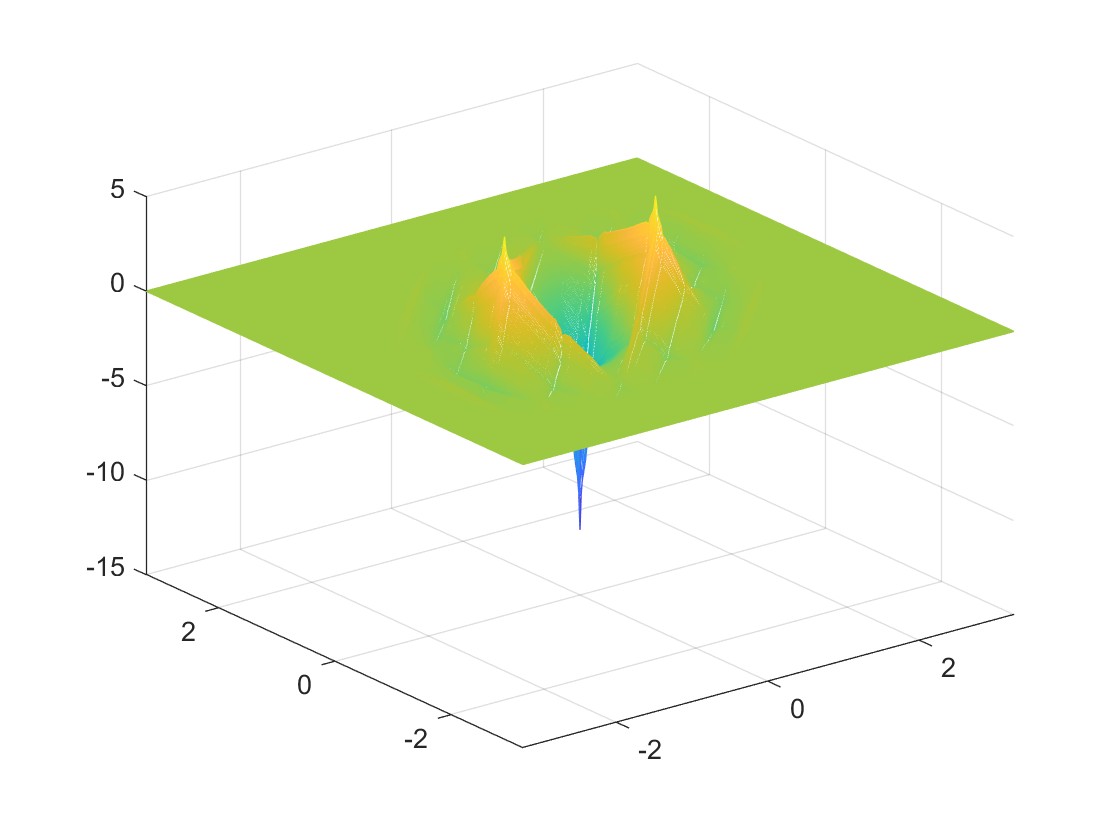} \caption{$\frac{\partial^2 \phi}{\partial x^2} \in C^0(\R^2)$}
	\end{subfigure}
	\begin{subfigure}[b]{0.22\textwidth} \includegraphics[width=\textwidth,height=0.8\textwidth]{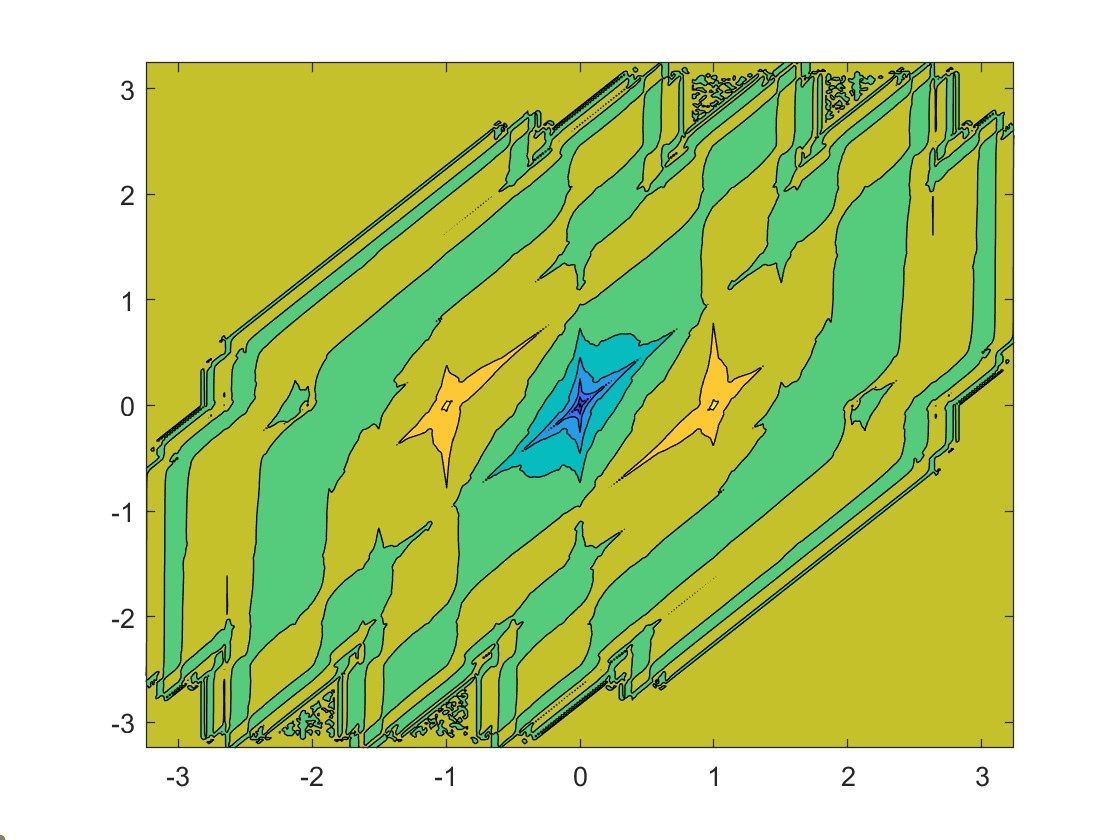}
		\caption{Contour of $\frac{\partial^2 \phi}{\partial x^2}$}
	\end{subfigure}\\
	\begin{subfigure}[b]{0.24\textwidth} \includegraphics[width=\textwidth,height=0.8\textwidth]{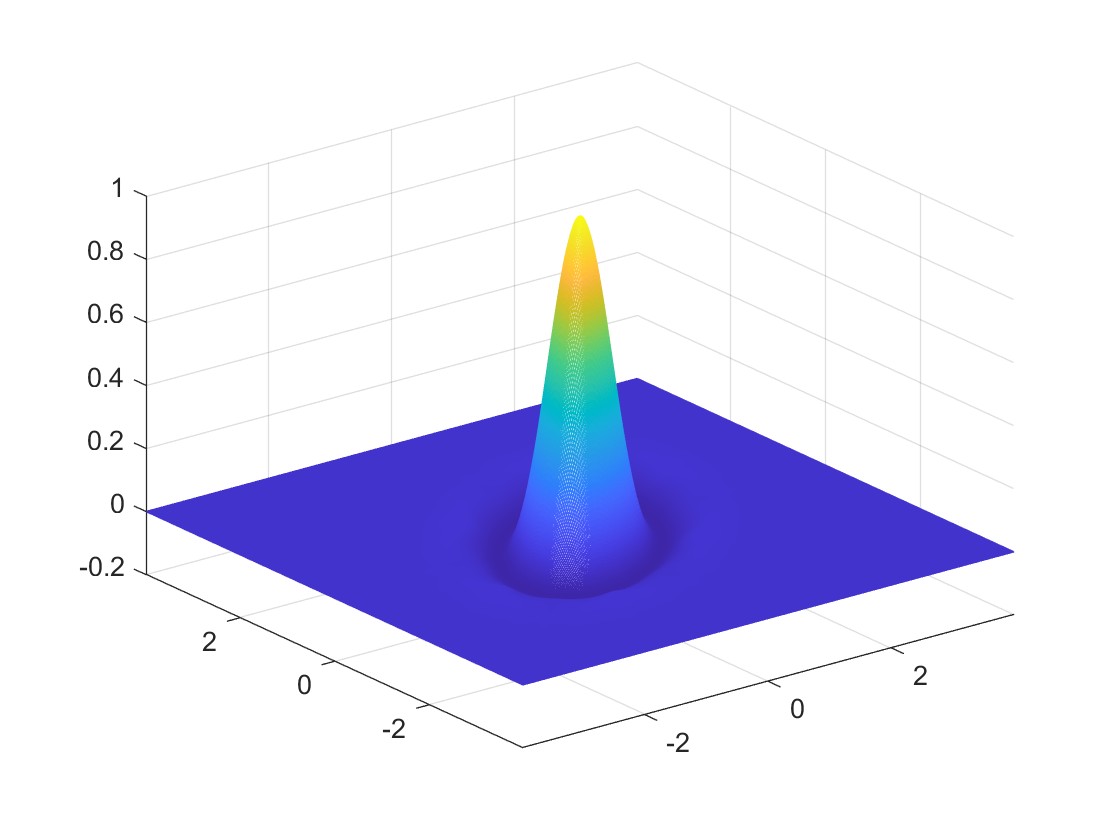}
		\caption{$\phi\in \mathscr{C}^2(\R^2)$}
	\end{subfigure}
	\begin{subfigure}[b]{0.22\textwidth} \includegraphics[width=\textwidth,height=0.8\textwidth]{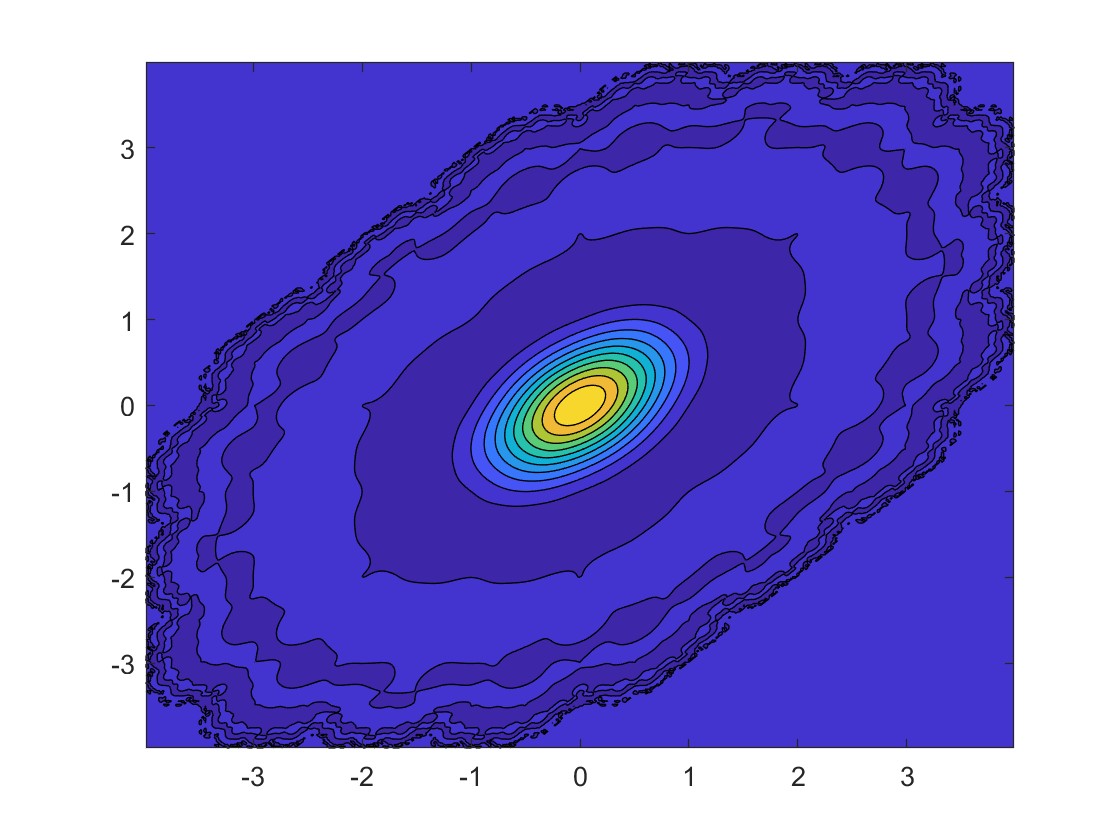}
		\caption{Contour of $\phi$}
	\end{subfigure}
	\begin{subfigure}[b]{0.24\textwidth} \includegraphics[width=\textwidth,height=0.8\textwidth]{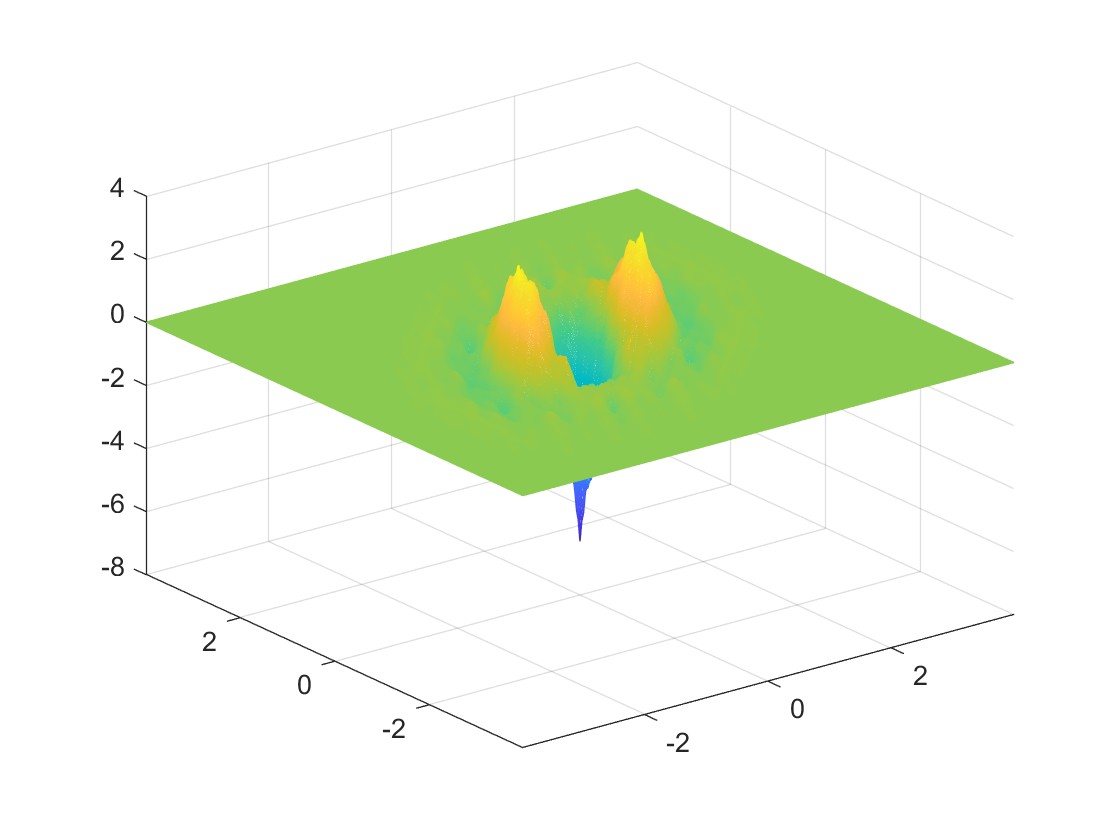} \caption{$\frac{\partial^2 \phi}{\partial x^2} \in C^0(\R^2)$}
	\end{subfigure}
	\begin{subfigure}[b]{0.22\textwidth} \includegraphics[width=\textwidth,height=0.8\textwidth]{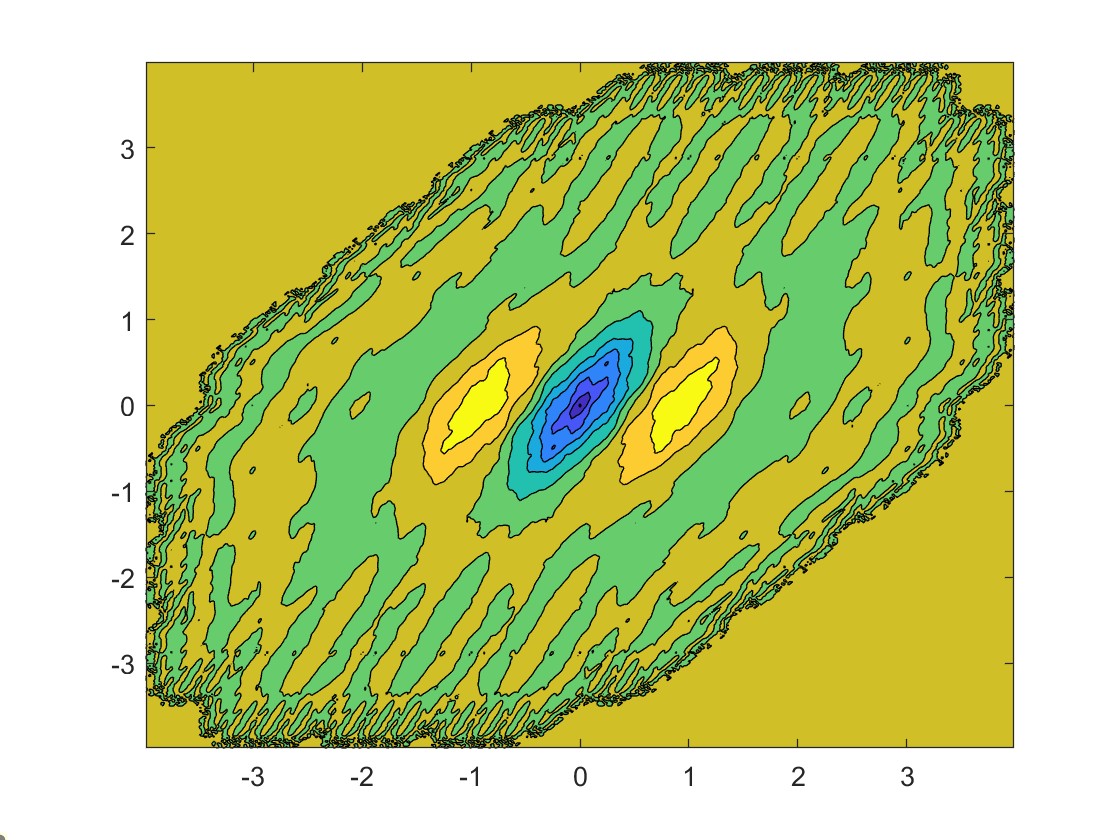}
		\caption{Contour of $\frac{\partial^2 \phi}{\partial x^2}$}
	\end{subfigure}
	
\caption{
		The first row is for the first choice in Example~\ref{ex2}:
		(A) is the graph of the interpolating $4I_2$-refinable function $\phi\in C^2(\R^2)$ and (B) is its contour.
		(C) is the graph of the partial derivative $\frac{\partial^2 \phi}{\partial x^2} \in C^0(\R^2)$, and (D) is the contour of $\frac{\partial^2 \phi}{\partial x^2}$.
		The second row is for the second choice in Example~\ref{ex2}:
		(E) is the graph of the interpolating $4I_2$-refinable function $\phi\in C^2(\R^2)$ and (F) is its contour.
		(G) is the graph of the partial derivative $\frac{\partial^2 \phi}{\partial x^2} \in C^0(\R^2)$, and (H) is the contour of $\frac{\partial^2 \phi}{\partial x^2}$.
	}\label{fig:ex2}
\end{figure}

\end{exmp}

%\textcolor{red}{you should add or delete the references to make it more balanced and complete.}

%\newpage

\end{document}